\documentclass[review,10pt]{elsarticle}
\usepackage[utf8]{inputenc}
\usepackage{microtype}
\usepackage[margin=2.5cm]{geometry}
\usepackage{setspace}
\usepackage[none]{hyphenat}
\usepackage[pagewise]{lineno}
\usepackage[english]{babel}
\usepackage{color}
\usepackage{booktabs}
\usepackage{longtable}
\usepackage[most]{tcolorbox} 
\usepackage{ragged2e}
\usepackage[overload]{empheq}
\usepackage[justification=justified]{caption}
\usepackage{enumitem}

\usepackage{amsmath,amsfonts,amssymb,mathtools}
\usepackage{cases}
\usepackage{multirow}
\usepackage{graphicx,float}
\usepackage{algpseudocode}
\usepackage{subcaption}
\usepackage{eurosym}
\usepackage{amstext} 
\usepackage{arydshln}
\usepackage{flushend}
\usepackage{threeparttable}
\usepackage{url}

\usepackage{float}

\usepackage[]{hyperref}
\hypersetup{colorlinks, linkcolor=blue, citecolor=blue, urlcolor=blue}

\usepackage[ruled,vlined,linesnumbered]{algorithm2e}
\DeclareRobustCommand{\officialeuro}{%
  \ifmmode\expandafter\text\fi
  {\fontencoding{U}\fontfamily{eurosym}\selectfont e}}
\usepackage{framed} 
\immediate\write18{%
makeindex Main_File.nlo -s nomencl.ist -o Main_File.nls 
}
\usepackage{multicol} 
\usepackage{nomencl}[noprefix] 
\makenomenclature
\renewcommand*\nompreamble{\begin{multicols}{2}}
\renewcommand*\nompostamble{\end{multicols}}

\def\tsc#1{\csdef{#1}{\textsc{\lowercase{#1}}\xspace}}
\tsc{WGM}
\tsc{QE}
\tsc{EP}
\tsc{PMS}
\tsc{BEC}
\tsc{DE}

\begin{document}

    \justifying
    \begin{frontmatter}   
        \title{Optimal Sizing of Community Photovoltaic and Battery Energy Storage Systems with Second-Life Batteries in Peer-to-Peer Energy Communities}

    
        \author[1,3]{Júlia Monar Pigem\corref{mycorrespondingauthor}}
        \cortext[mycorrespondingauthor]{Corresponding author}
        \ead{julia.monar@estudiantat.upc.edu}
        \author[1,2]{Fernando García-Muñoz}
        \author[1,2]{Natalia Jorquera Bravo}
        \author[1]{Joaquín Aballay Araya}
        \author[4]{Vicente Castro Burgos}

        \address[1]{University of Santiago of Chile (USACH), Faculty of Engineering, Industrial Engineering Department, Chile}
        \address[2]{University of Santiago of Chile (USACH), Faculty of Engineering, Program for the Development of Sustainable Production Systems (PDSPS), Chile}
        \address[3]{Universitat Politècnica de Catalunya (UPC), Escola Tècnica Superior d'Engyinyeria Industrial de Barcelona, Spain}
        \address[4]{University of Santiago of Chile (USACH), Faculty of Engineering, Electrical Engineering Department, Chile}

    \begin{abstract}
    This article presents a mixed-integer second-order cone programming model to determine the optimal sizing of a community-shared photovoltaic and battery energy storage system (PV-BESS) within a peer-to-peer (P2P) energy trading framework. The model accounts for heterogeneous users who may already own individual PV or PV-BESS systems and aims to enhance the overall energy autonomy of the energy community. A key feature of the model is the explicit comparison between first-life (FL) and second-life (SL) battery technologies, incorporating their respective degradation dynamics into investment and operational decisions, and the technical feasibility by considering constraints of a low-voltage distribution network. The proposed formulation is tested on the reduced equivalent of the IEEE European low-voltage network. Results show that the most influential factors in the adoption of a shared BESS are: (i) the market cost of battery technologies, (ii) electricity tariffs, particularly purchase prices, and (iii) the degradation characteristics of the chosen technology. Secondary factors, such as DER penetration among users and the community’s peak demand, have a lesser impact. The analysis further suggests that SL batteries could become a cost-effective alternative to FL technologies if their degradation performance improves or their capital cost is significantly reduced.
    \end{abstract}


\end{frontmatter}

\section{Introduction}\label{Intro_Chap}

The increasing competitiveness of distributed energy resources (DERs), such as photovoltaic (PV) systems and battery energy storage systems (BESS), has accelerated their adoption among residential users, prompting the electrical system to adapt accordingly~\cite{2}. As more households integrate these technologies, energy communities have emerged as a promising framework to enhance local energy autonomy, improve system resilience, and reduce reliance on the traditional power grid~\cite{bollen2011integration}. In this regard, peer-to-peer (P2P) energy trading has gained attention as an effective mechanism for local energy exchange and better utilization~\cite{11}. P2P markets allow prosumers to directly sell or share their excess generation with neighboring users, promoting efficient utilization of renewable resources~\cite{Wang2018}. However, this new operational paradigm modifies traditional consumption patterns, as generation peaks often misalign with demand~\cite{Barata2025}. In this context, BESS becomes critical in storing surplus PV energy for later use, enabling the full potential of P2P exchanges~\cite{Tejero2024G}. Thus, BESS play a key role in reducing grid dependency and enabling responsive energy management~\cite{37}. Their ability to shift energy usage across time enhances self-consumption and facilitates adaptation to time-varying generation and price signals~\cite{Abdullah2021}. This makes them essential components for modern distributed power systems based on DER coordination and market-based energy exchanges.

In addition to user-owned DERs, some energy communities aim to invest in shared assets such as community-scale PV systems and communal BESS~\cite{Goulart2024}. These installations provide collective benefits through economies of scale and broaden access to renewable generation and storage capacity~\cite{Rowe2023}. Their planning and operation, however, must account for the particularities of community dynamics and new trading mechanisms like P2P. In particular, communal BESS can enhance the community’s ability to balance generation and demand across users, thereby facilitating more efficient and equitable energy exchanges. Nevertheless, this added flexibility comes at a cost; BESS still represents a significant capital investment, especially at the community scale~\cite{Zurita2018}. In this context, repurposed batteries from electric vehicles (EV), known as second-life (SL) batteries for stationary applications, have emerged as a cost-effective alternative~\cite{Demirci2023}. Although more affordable, they suffer from reduced performance and accelerated degradation~\cite{Gao2023}. This trade-off between investment cost and long-term usability raises key questions about their suitability for community applications, and highlights the need for planning tools that can evaluate their contribution to autonomy and cost-efficiency~\cite{Iqbal2023}. 

In this context, evaluating the feasibility and performance of community-scale PV-BESS systems becomes essential, particularly in configurations where users already own individual assets and participate in local energy exchanges. Given the economic challenges BESS investments pose, especially under dynamic electricity market conditions, exploring alternative technologies that can reduce upfront costs while maintaining acceptable performance is crucial. SL batteries offer a promising yet uncertain option due to their lower cost but potentially limited lifespan. Accordingly, optimization-based planning tools that assess the optimal sizing of shared PV and BESS systems, accounting for different battery chemistries, degradation patterns, and usage contexts, are particularly relevant. Such approaches enable a systematic evaluation of the trade-offs between investment cost and long-term performance, especially when comparing first- and second-life technologies for collective use in energy communities. In this work, such a model is developed to quantify these trade-offs between both technologies when deployed in a collective setting, providing insights into the conditions under which SL solutions may become viable for community applications. Unlike previous studies, it explicitly addresses the integration of P2P energy trading, realistic battery degradation dynamics, and network constraints, thereby filling a critical gap in the planning of community-scale PV-BESS systems. 

The remainder of the paper is organized as follows. Section~\ref{sec:literature} reviews the relevant literature. Section~\ref{sec:Methodology} describes the methodology used to construct the demand, generation, and degradation scenarios. Section~\ref{sec:Opt_Model} introduces the proposed optimization model. Section~\ref{sec:results} presents the case study and discusses the computational results. Finally, Section~\ref{sec:conclusions} summarizes the main findings and outlines directions for future research.

\section{Literature Review}\label{sec:literature}
Several studies have explored the optimal sizing of PV-BESS systems in grid-connected contexts, focusing on economic feasibility, self-consumption enhancement, and operational efficiency. For instance, the study in~\cite{Zhou2025} proposes a mixed-integer linear programming (MILP) model to optimize PV and BESS configurations in a university campus, demonstrating their potential to increase electricity autonomy. Similarly, the works in~\cite{Alic2024} and~\cite{Aripriharta2025} address PV-BESS design for residential systems that include EVs, employing Particle Swarm Optimization (PSO) and Grey Wolf Optimizer (GWO) algorithms to solve the problem across multiple urban contexts. These studies arrive at comparable conclusions, reinforcing the potential of PV-BESS systems to improve local energy resilience and reduce grid dependency. This is further supported by~\cite{Linssen2017} and~\cite{Tsai2020}, who show that PV-BESS integration enhances self-consumption and financial performance. Additionally, the work in~\cite{Mahmud2022} uses HOMER software to identify cost-optimal microgrid configurations in a Bangladeshi community, while~\cite{Venkatasubramanian2025} develops a framework for capacity expansion in active microgrids. Demand-side flexibility has also been incorporated in~\cite{Li2022} to maximize renewable energy utilization, while~\cite{Abdullah2021} analyzes the economic viability of BESS under different sizing and discharge profiles, highlighting the impact of battery capacity on cost-effectiveness. Although these contributions demonstrate that PV-BESS systems can enhance self-sufficiency and resilience across different settings, several key aspects remain unaddressed in previous works. In particular: (i) the potential role of emerging local electricity markets, such as P2P trading schemes, which rely heavily on storage for their operation; (ii) the long-term impact of battery degradation; and (iii) the diversity of available battery technologies, including comparisons between first-life (FL) and SL batteries, are all could be critical for reliable and cost-effective system planning.

In the above context, some researchers have begun exploring the role of emerging local electricity markets as a fundamental element in designing and operating a shared PV-BESS system. For instance,\cite{Aoun2024} compares blockchain-based P2P energy trading with traditional incentive mechanisms such as Net Energy Metering and Feed-in Tariffs in rural microgrids, showing that P2P schemes enable more flexible and economically attractive real-time energy exchanges. Similarly,\cite{Bokopane2024} develops a model for a grid-connected EV charging station microgrid that integrates renewable energy, battery storage, and P2P energy sharing, demonstrating that dynamic pricing mechanisms can effectively facilitate prosumer participation while reducing operational costs. Building on this, \cite{Kachhad2024} proposes a techno-economic framework for the design of Renewable Energy Communities, emphasizing the role of technical facilitators and component sizing in enhancing financial outcomes through P2P trading. The study in \cite{Zhang2024} presents an optimization model for BESS sizing under different electricity market designs, including P2P and energy storage sharing, and evaluates alternative ownership and utilization strategies. Similarly,~\cite{Jose2024} explores a residential microgrid where prosumer selection is optimized via a genetic algorithm, finding that a balanced composition of consumers and prosumers improves the overall economic performance of the system. However, while these studies incorporate market-related mechanisms and highlight the value of local energy trading, they still overlook critical aspects such as battery degradation and the diversity of available BESS technologies.

As a result of the above gap, the research community has begun to evaluate different battery technologies and the effects of degradation in the context of shared PV-BESS systems. For instance, \cite{Torkashvand2020} analyzes the life cycle costs of lithium-ion and lead-acid batteries for islanded microgrids, modeling degradation by considering the impact of charge–discharge cycles and operational stress on capacity over time. Similarly, \cite{Stamantis2019} compares lithium-ion, lead-acid, sodium-sulfur, and vanadium redox flow batteries for PV systems aimed at enhancing resilience to grid outages, evaluating them based on efficiency, cost, depth of discharge, and lifetime. In~\cite{Bandyopadhyay2020}, degradation is incorporated as capacity loss over time when optimizing system sizing. The study in~\cite{Salman2021} applies depth-of-discharge modeling to extend battery lifetime and reduce costs in a wind-BESS microgrid. Lastly,~\cite{Mohamed2022} develops a tool for techno-economic analysis of small-scale PV-BESS systems that includes battery degradation, concluding that accounting for degradation renders such systems economically unviable under current commercial prices, while highlighting the critical role of degradation in the decision-making process.

This conclusion is further supported by several studies that have highlighted the technical and economic limitations of FL batteries in stationary applications, often emphasizing the presence of significant cost-related barriers. For instance,\cite{Barata2025} finds that standalone community-scale BESS are not yet economically viable at current technology and price levels unless complemented by additional revenue streams. Similarly, \cite{Goulart2024} models the optimal sizing and operation of BESS while quantifying potential revenues from ancillary services such as frequency regulation, voltage support, energy time-shifting, and backup power, showing that economic feasibility strongly depends on future declines in battery prices. In~\cite{Zurita2018}, a hybrid renewable power plant is evaluated, concluding that profitability requires a 60–90\% reduction in storage costs. Along the same lines,~\cite{Lopez-Lorente2021} analyzes centralized grid-scale BESS supporting distribution networks with high PV penetration and finds that current market conditions remain insufficient to justify investment under existing regulatory and pricing frameworks.

In light of these economic limitations, SL batteries, repurposed from EVs for stationary applications, have emerged as a potential alternative to reduce investment costs and improve the economic viability of BESS-based systems. In this regard, \cite{Demirci2023} optimizes the sizing and techno-economic performance of 100\% renewable energy systems using SL batteries, explicitly modeling degradation through capacity fade and cycle life. The study evaluates lifecycle costs and shows that SL batteries can reduce expenses sustainably while maintaining reliable storage performance. Similarly, \cite{Silvestri2021} demonstrates the technical feasibility and economic benefits of implementing a SL BESS in a ceramic manufacturing plant, enabling peak shaving, lowering electricity demand charges, and increasing solar self-consumption. The comparison between first-life and SL EV batteries is further explored in~\cite{Bai2024}, which presents a techno-economic assessment of isolated microgrids using a reliability-oriented iterative design framework. The model includes a semi-empirical battery degradation representation to ensure demand satisfaction and realistic lifetime estimates. Along the same lines,~\cite{Wangsupphaphol2023} evaluates a hybrid system for home EV charging based on solar PV and both FL and SL batteries, using a detailed degradation model accounting for cycle count, depth of discharge, temperature, and state of charge. The results show strong economic returns over a 10-year horizon, positioning SL batteries as a cost-effective and sustainable alternative. In~\cite{Terkes2023}, hybrid power systems combining FL and SL batteries are analyzed, indicating that while SL batteries offer significant cost and environmental benefits, larger storage capacities may be required to match the performance of new units. Additional evidence of the benefits of SL batteries in community-scale applications is presented by~\cite{Deng2022}, who report that repurposed LiFePO$_4$ batteries can reduce capital costs by up to 94\% compared to FL lithium-ion units, while still providing effective peak shaving and grid support in residential microgrids. These findings highlight the potential of SL batteries to enhance the viability of energy projects, particularly in low-income or cooperative-based community energy systems.

The previous literature reveals a growing interest in optimization-based studies that compare FL and SL batteries across various stationary applications, including residential~\cite{Wangsupphaphol2023}, industrial~\cite{Silvestri2021}, and community-scale systems~\cite{Deng2022}. This is explained mainly because FL batteries, despite having decreased pricing over the last years, still have a remaining gap in competition against electricity market prices; hence, they have promoted the SL battery studies. Likewise, these studies increasingly acknowledge the importance of explicitly modeling battery degradation, particularly relevant in the case of SL batteries, which, although more affordable in terms of investment cost, exhibit accelerated degradation. However, despite these advances, there is still a lack of optimization models specifically designed for community-based energy systems that integrate emerging market structures such as P2P trading, while simultaneously comparing FL and SL battery technologies under realistic degradation assumptions. Moreover, ensuring the technical feasibility of such systems requires the inclusion of network constraints, which are often neglected in prior works, even when dimensioning shared community assets. In addition, existing studies that evaluate battery degradation typically rely on predefined or average usage patterns, rather than on degradation profiles derived from optimized consumption schedules. Accordingly, this article aims to address these gaps by proposing the following contributions:

\begin{itemize}
    \item A mixed-integer second-order cone programming (MISOCP) model to determine the optimal sizing and technology selection of a community-scale PV-BESS system within a P2P energy trading framework. The model considers users owning individual DERs and incorporates low-voltage distribution network constraints.

    \item The model explicitly compares FL and SL battery technologies by integrating their respective degradation dynamics into the sizing and operational decisions. To do so, a separate optimization-based framework is employed to simulate typical daily BESS usage and estimate degradation profiles (SoC profile), which are fed into the main MISOCP formulation.

    \item A sensitivity analysis was conducted to explore the economic competitiveness of FL and SL batteries across varying market contexts. By adjusting electricity prices and technology costs, the analysis identifies the conditions under which each battery type becomes economically viable and evaluates when degradation effects outweigh capital cost advantages.
\end{itemize}

\section{Methodology} \label{sec:Methodology}
\begin{figure}[H]
    \centering
    \includegraphics[width=0.8\textwidth]{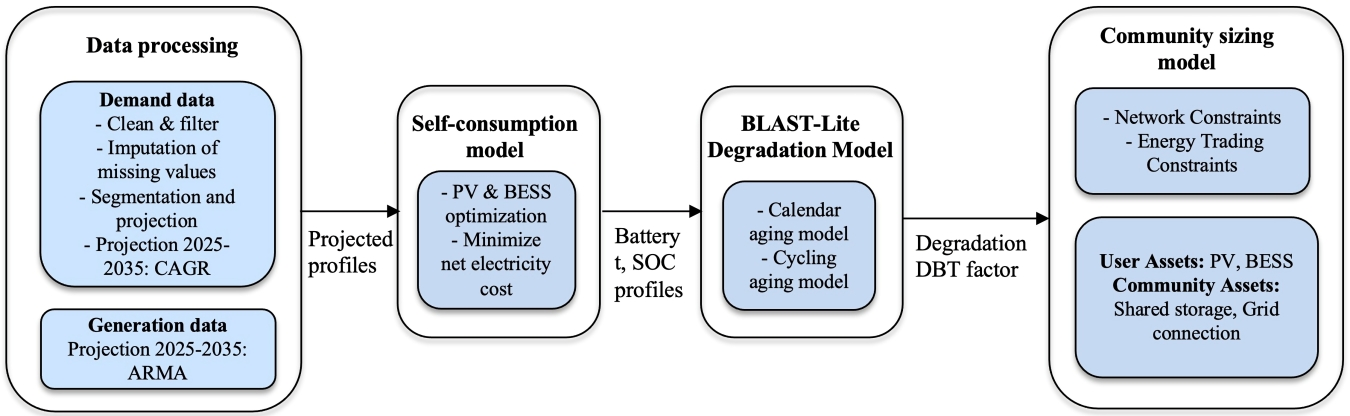} 
    \caption{Methodology map}
    \label{fig:mapa_methodology}
\end{figure}

This section introduces the methodological framework developed to support the optimal sizing of community-scale PV-BESS systems under detailed technical and economic considerations. The overall process is summarized in Figure~\ref{fig:mapa_methodology}, illustrating how different modeling layers are integrated to generate the necessary input data for the final MISOCP-based sizing model.

The process begins with the preparation of residential electricity demand and solar irradiance profiles, which are cleaned, segmented, and projected for the 2025–2035 horizon. These profiles are then used as input for a user-level self-consumption model that optimizes PV and BESS operation for a single user by minimizing net electricity costs. The output of this model, specifically the battery’s hourly state-of-charge (SoC) profile, is used to estimate degradation through the BLAST-Lite simulator, a MATLAB-based tool that evaluates the degradation of the different BESS technologies considered in this study. Finally, the resulting degradation indicators, together with the projected demand and generation profiles, are integrated into a community-level MISOCP model.


\subsection{Data processing}

The dataset used in this study~\cite{garcia2024exploring} contains hourly electricity consumption profiles from both residential and commercial users in a rural Spanish distribution network. Subsequently, only those residential profiles with an average annual consumption above a defined threshold were retained, following the criteria proposed in~\cite{VanderPlas2016}. Given the presence of missing data, common in real-world metering systems, gaps were imputed using the k-nearest neighbors method~\cite{fadlil2022knn}, and redundant entries were removed. The remaining profiles were then classified into constant and variable consumption types. For the latter, long-term demand was projected over a 10-year horizon using a seasonal cyclic approach~\cite{Taylor2010}. This process resulted in 13 representative residential profiles, each consisting of 8760 hourly values, from which a subset aligned with a predefined reference profile was selected for the optimization model. Next step has been to estimate projections of electricity demand over a 10-year horizon. To achieve this, the results obtained in~\cite{jimenez2023} were used, as the reference study provides projections of demand for a rural community archetype in Spain that aligns with the consumption profiles contained in our data set. Plus, it accounts for the rise in energy consumption driven by the growing adoption of electric vehicles and heat pumps.


Based on these projections and the optimization model developed in this work, annual estimates of electricity demand were established for the years 2025, 2030, and 2040. To quantify the expected rate of growth between these periods, the Compound Annual Growth Rate (CAGR) was calculated. This approach is widely recognized for its ability to consistently characterize long-term growth patterns over defined time intervals, as highlighted in~\cite{chen2023}. The CAGR was determined for each interval using its standard mathematical formulation. To enhance the robustness of the analysis, given the 10-year temporal framework adopted in this study, a constant annual growth rate was assumed. This rate corresponds to the arithmetic mean of the CAGR values calculated for each time interval and it provides the most coherent and stable representation of the projected demand evolution. This allows the projection of energy consumption over the next 10 years. In addition, normalized consumption profiles are calculated to compare the relative consumption of each household~\cite{zhang2003}, and a visualization is generated with the results of the total demand projections. 

To model and project solar generation, the AutoRegressive Moving Average (ARMA) method is employed. This statistical technique is well-suited for time series forecasting, particularly when the data exhibits temporal autocorrelation~\cite{Box2015}. The ARMA model combines two components: the autoregressive (AR) part, which captures the influence of past values, and the moving average (MA) part, which accounts for past forecast errors. Historical hourly solar generation data is used to fit the ARMA parameters, ensuring the model captures both the seasonal patterns and random fluctuations present in the data~\cite{Singh2019ARMA}. Once trained, the ARMA model generates synthetic time series that preserve the statistical properties of the original dataset, providing a realistic basis for scenario analysis and simulation within the broader energy system modeling framework. For 10 year projections, the first year of historical data is used as the input to the ARMA model, which is then iteratively simulated 10 times to generate the time series.

\subsection{Self consumption model}

The estimation of long-term degradation for different BESS technologies begins with a simplified optimization model that simulates the energy management of a single user equipped with PV generation and a BESS. The model aims to determine the optimal operation profile of BESS operation over a one-year horizon at hourly resolution (8760 time steps), assuming perfect foresight of demand and generation. This operation profile is subsequently used as input to a degradation analysis with the BLAST-Lite model to compute the degradation factor \( DBT \) associated with each battery technology, which is then incorporated into the Community PV-BESS Sizing Model. Thus, in this first stage, the objective function in~\eqref{Single_Op_Model} minimizes the user's net electricity cost over the horizon, computed as the total cost of energy purchased from the grid \( p^{bg}_{t} \) at price \( \lambda^{bg}_{t} \), minus the revenue from selling surplus energy \( p^{sg}_{t} \) at price \( \lambda^{sg}_{t} \).
\begin{align}
        &\min \enspace z = \sum_{t \in \mathcal{T}} \left( \lambda^{bg}_{t} \, p^{bg}_{t} - \lambda^{sg}_{t} \, p^{sg}_{t} \right)\label{Single_Op_Model}
\end{align}

\begin{subequations}
\begin{align}
    &pv_{t} - PL_{t} + ds_{t} - ch_{t} + p^{bg}_{t} - p^{sg}_{t} = 0 && \forall t \in \mathcal{T} \label{SC_Model_Balance} \\
    &pv_{t} = \Gamma^{pv} \cdot PG^{max}_{t} && \forall t \in \mathcal{T}  \label{SC_Model_PV} \\
    &soc_{t} = soc_{t-1} + \left( \varphi^{ch} \, ch_{t} - \frac{1}{\varphi^{ds}} \, ds_{t} \right) \Delta t && \forall t \in \mathcal{T} \label{SC_Model_SOC}  \\
    &SOC^{min} \cdot \Gamma^{bt} \leq soc_{t} \leq SOC^{max} \cdot \Gamma^{bt} && \forall t \in \mathcal{T}  \label{SC_Model_SOC_limit} \\
    &ch_{t} \leq \frac{1}{SB} PB \cdot w_{t} && \forall t \in \mathcal{T}  \label{SC_Model_ch} \\
    &ds_{t} \leq \frac{1}{SB} PB \cdot (1 - w_{t}) && \forall t \in \mathcal{T}  \label{SC_Model_ds} 
\end{align}
\end{subequations}

The power balance constraint~\eqref{SC_Model_Balance} ensures that, at each time step \( t \in \mathcal{T} \), the user's electricity demand \( PL_{t} \) is met by PV generation \( pv_{t} \), battery discharging \( ds_{t} \), and grid imports \( p^{bg}_{t} \), and is offset by battery charging \( ch_{t} \) and grid exports \( p^{sg}_{t} \). PV generation is modeled in~\eqref{SC_Model_PV} as a known installed capacity \( \Gamma^{pv} \) scaled by a normalized irradiance profile \( PG^{max}_{t} \). The battery’s SoC dynamics are captured by~\eqref{SC_Model_SOC}, which accounts for charging and discharging efficiencies \( \varphi^{ch} \) and \( \varphi^{ds} \). The SoC is bounded within minimum and maximum thresholds, scaled by the installed battery capacity \( \Gamma^{bt} \), as specified in~\eqref{SC_Model_SOC_limit}. The operational window is restricted to 10–90\% in order to avoid the excessive degradation observed at very low and very high SoC levels~\cite{Vermeer}.  Constraints~\eqref{SC_Model_ch} and~\eqref{SC_Model_ds} use the binary variable \( w_t \) to enforce mutually exclusive charging and discharging modes, each limited by the rated battery power \( PB \).

\subsection{Matlab BLAST-Lite Model}

Battery degradation is estimated using BLAST-Lite, a simulation tool to evaluate the lifetime performance of lithium-ion batteries in both stationary and mobile applications~\cite{blastlite}. The model takes as input the SOC profiles generated by the self-consumption model, along with the corresponding time and temperature data. It incorporates four predefined battery types, which are summarized in Table~\ref{tab:BLAST modelos de bateria}, each with specific chemistry and design characteristics. The cathode materials include lithium iron phosphate (LFP), lithium manganese oxide (LMO), nickel cobalt aluminum oxide (NCA), and nickel manganese cobalt oxide (NMC), while the anode materials are graphite (Gr), graphite with silicon additive (GrSi), and lithium titanate (LTO). An LMO-based SL battery was selected for the degradation study, both because this chemistry is commonly used in SL applications~\cite{Sarker} and due to the availability of resources.

\begin{table}[h!]
\centering
\caption{Battery types included in BLAST LITE Matlab Model}
\begin{tabular}{|c|c|c|c|c|c|p{8cm}|}
\hline
\textbf{Type} & \textbf{Cathode} & \textbf{Anode} & \textbf{Battery life stage} & \textbf{Capacity [Ah]} &\textbf{Reference}\\
\hline
LFP/Gr & LFP & Graphite & FL & 250 & \cite{Gasper2023}\\
LMO/Gr & LMO & Graphite & SL & 66 & \cite{Braco2020}\\
NMC/Gr & NMC & Graphite & FL & 75 &\cite{Gasper2023}\\
NMC/LTO & NMC & Lithium titanate & FL & 10 &\cite{nmc/lto}\\
\hline
\end{tabular}
\label{tab:BLAST modelos de bateria}
\end{table}

  The calculation of battery degradation varies slightly between cell types due to differences in chemistry, design, and operational sensitivity~\cite{Gasper2023}. As a result, models must be specifically tailored to each battery. However, all models rely on a common set of stressors to quantify degradation, including time, cycling activity, depth of discharge (DOD), charge/discharge rate (C-rate), temperature, SOC, and anode potential. In cases with dynamic cycling, the Rainflow counting algorithm is applied to differentiate full from partial cycles and to quantify fatigue-related damage.

  Most battery types in BLAST-Lite, such as LFP/Gr or LMO/Gr follow a semi-empirical degradation model that separates capacity loss into calendar and cycling components. Calendar aging is modeled using a power-law function of time, with degradation rates defined through simplified Arrhenius and Tafel relationships. These expressions capture the effects of temperature and anode potential ($U_a$). Cycling degradation is modeled as a power-law function of Full Equivalent Cycles (FEC), with rates accelerating under high or low temperatures, increased DOD, and higher C-rates. The equations are presented in the appendix, and all model parameters are obtained by fitting them to experimental datasets from each specific battery type. In the case of SL batteries, LMO/Gr, an additional equation is incorporated into the degradation model to constrain the remaining capacity, accounting for its lower initial SOH. 

  The model for NMC/Gr follows a similar structure, incorporating exponential dependencies on temperature and $U_a$ for calendar aging, alongside non-linear dependencies on DOD, temperature, and C-rate for cycling degradation. Similarly, the degradation model for NMC/GrSi accounts for capacity fade driven by both calendar and cycling mechanisms. Calendar aging is described by a combination of exponential and polynomial dependencies on normalized temperature, $U_a$ and SOC while cycle aging is modeled using exponential terms related to discharge C-rate and temperature. In both models, the degradation contributions are aggregated using power-law functions of time and equivalent EFCs. All relevant parameters and detailed equations are provided in the appendix.

  In contrast, the NMC/LTO model introduces both irreversible capacity loss and reversible capacity gain during calendar aging. These two components are modeled separately as functions of temperature and SOC using exponential terms. Cycling degradation is treated independently, depending on DOD, C-rate and temperature. All three contributions, calendar loss, calendar gain, and cycling loss, are accumulated over time or EFCs using power-law relationships. The corresponding equations are listed in the appendix.

\section{Community PV-BESS Sizing Model}\label{sec:Opt_Model}
This section presents the MISOCP model, where the main notation used in this paper is summarized below () for quick reference, while other symbols are defined as needed throughout the text.

\begin{table}[ht]   
\small
\begin{framed}
\nomenclature[24]{\textbf{Sets}}{}
\nomenclature[25]{$i \in \Omega$}{Index of buses.}
\nomenclature[26]{$(i,j) \in \mathcal{L}$}{Set of distribution lines.}
\nomenclature[27]{$i \in \Omega_A$}{Subset of buses with active users, such that $\Omega_A \subseteq \Omega$}
\nomenclature[28]{$i \in \Omega_S$}{Subset of buses with substation connection, such that $\Omega_S \subseteq \Omega$}
\nomenclature[29]{$t \in \mathcal{T}$}{Set of time steps (e.g., hours).}
\nomenclature[30]{$r \in \mathcal{R}$}{Set of years in the planning horizon.}
\nomenclature[31]{$f \in \mathcal{F}$}{Set of battery technologies.}

\nomenclature[32]{\textbf{Parameters}}{}
\nomenclature[33]{$R_{i,j}$}{Resistance of line $(i,j)$ [$\Omega$].}
\nomenclature[34]{$X_{i,j}$}{Reactance of line $(i,j)$ [$\Omega$].}
\nomenclature[36]{$PB$}{Maximum power (kW) of user BESS.}
\nomenclature[37]{$PB^{cm}$}{Maximum power (kW) of community BESS.}
\nomenclature[38]{$PSG_i^{\max}$}{Max active power sold to the grid.}
\nomenclature[38]{$PBG_i^{\max}$}{Max active power bought from the grid.}
\nomenclature[39]{$PL_{i,t,r}$}{Active power Load of user $i$ at time $t$, year $r$ [kW].}
\nomenclature[39]{$QL_{i,t,r}$}{Reactive power Load of user $i$ at time $t$, year $r$ [kvar].}
\nomenclature[40]{$I^{pv}$}{Marginal cost of PV-shared system [€/kW].}
\nomenclature[41]{$I^b_f$}{Marginal cost of BESS-shared technology $f$ [€/kWh].}
\nomenclature[42]{$SOC^{\min}$}{Min and max SoC bounds [\%].}
\nomenclature[42]{$SOC^{\max}$}{Max and max SoC bounds [\%].}
\nomenclature[43]{$DBT_{r,f}$}{Degradation factor for BESS tech $f$ in year $r$.}
\nomenclature[44]{$\eta$}{Charge/discharge efficiencies.}
\nomenclature[45]{$BT_i$}{Binary, 1 if user $i$ has BESS.}
\nomenclature[46]{$\text{SB}$}{Base power used for per-unit normalization.}

\nomenclature[49]{\textbf{Variables}}{}
\nomenclature[50]{$\gamma_i^{pv}$}{Installed PV capacity at node $i$ [kW].}
\nomenclature[51]{$\gamma_{i,f}^{bt}$}{Installed BESS capacity of type $f$ at user $i$ [kWh].}
\nomenclature[53]{$pv_{i,t,r}$}{PV generation at node $i$, time $t$, year $r$ [kW].}
\nomenclature[54]{$p_{i,t,r}^{bg}$}{Power bought from grid [kW].}
\nomenclature[54]{$p_{i,t,r}^{sg}$}{Power sold to grid [kW].}
\nomenclature[54]{$p_{i,t,r}^{sm}$}{Power sold to the local market [kW].}
\nomenclature[54]{$p_{i,t,r}^{bm}$}{Power bought from the local market [kW].}
\nomenclature[55]{$ch_{i,t,r}$}{Charging power [kW].}
\nomenclature[55]{$ds_{i,t,r}$}{Discharging power [kW].}
\nomenclature[56]{$soc_{i,t,r}$}{State of charge of BESS at node $i$ [
\nomenclature[57]{$w_{i,t,r}$}{Binary variable: 1 if charging at $i$ at $t,r$; 0 otherwise.}
\nomenclature[58]{$pv_{i,t,r}^{cm}$}{Shared PV power used at $i$ [kW].}
\nomenclature[59]{$ch_{i,t,r,f}^{cm}$}{Charging of community BESS.}
\nomenclature[59]{$ds_{i,t,r,f}^{cm}$}{Discharging of community BESS.}
\nomenclature[60]{$soc_{i,t,r,f}^{cm}$}{SoC of community BESS tech $f$ at $i$ at $t,r$.}
\nomenclature[61]{$y_{i,t,r}$}{Binary: 1 if selling at time $t,r$, 0 otherwise.}
\nomenclature[62]{$\Delta p_{i,t,r}$}{Net power injection at $i$ and $t,r$ [kW].}
\nomenclature[63]{$\Delta p_{i,t,r}^+$}{Surplus components of net injection.}
\nomenclature[63]{$\Delta p_{i,t,r}^-$}{Deficit components of net injection.}
\nomenclature[64]{$p\kappa_{i,t,r}^{sg}$}{Energy exported to neighbor $j$.}
\nomenclature[65]{$p\kappa_{i,t,r}^{bg}$}{Energy imported from neighbor $j$.}
\nomenclature[66]{$v_{i,t,r}$}{Squared voltage magnitude at node $i$ at $t,r$ [$\text{p.u.}^2$].}
\nomenclature[67]{$\ell^{i,j}_{t,r}$}{Squared current magnitude on line $(i,j)$ at $t,r$.}
\nomenclature[68]{$p^{i,j}_{t,r}$}{Active power flow on line $(i,j)$.}
\nomenclature[68]{$q^{i,j}_{t,r}$}{Reactive power flow on line $(i,j)$.}
\nomenclature[69]{$qg_{i,t,r}$}{Reactive power injection at node $i$.}

\printnomenclature
\end{framed}
\end{table}

The proposed optimization model considers a low-voltage radial DN representing an EC, which is modeled as a directed graph \( \mathcal{G} = (\mathcal{B}, \mathcal{L}) \), where \( \mathcal{B} \) denotes the set of buses and \( \mathcal{L} \) the set of distribution lines. A subset of buses \( \mathcal{A} \subseteq \mathcal{B} \) corresponds to active community members, some of whom are currently equipped with individual PV systems and/or BESS. These users already engage in P2P energy trading, aiming to increase local energy autonomy and reduce reliance on the upstream grid. To further support daily operations and enhance the community's resilience, the EC seeks to invest in a shared PV-BESS system to be installed at the point of common coupling \( S \in \mathcal{B} \), which connects the community to the upstream DN. The objective of the model is to determine the optimal capacity of the community PV-BESS system to be deployed, such that it supports the operation of the EC over a 10-year planning horizon \( \mathcal{Y} \), considering hourly resolution \( t \in \mathcal{T} \) for each day of the year. The optimization model is formulated as a Mixed-Integer Second-Order Cone Programming (MISOCP) problem, combining binary variables for discrete investment and operational decisions with convex conic constraints derived from the second-order relaxation of AC power flow equations.

  The investment cost includes a marginal cost \( I^{pv} \) per kW for the PV system and a technology-dependent marginal cost \( I_f^{bt} \) for each BESS option \( f \in \mathcal{F} \), allowing the model to select among different battery technologies with distinct cost-performance characteristics. This flexibility is particularly relevant in scenarios involving second-life batteries, where trade-offs between cost and degradation must be evaluated. The installed capacities are represented by \( \gamma^{pv}_i \) for PV and \( \gamma^{bt}_{i,f} \) for each BESS type. During operation, the community may exchange energy with the external grid. At each time step \( t \in \mathcal{T} \) and year \( r \in \mathcal{Y} \), energy can be purchased from the grid at a price \( \lambda_{t,r}^{bg} \) and sold at a price \( \lambda_{t,r}^{sg} \), with the respective exchange variables \( p\kappa^{bg}_{i,t,r} \) and \( p\kappa^{sg}_{i,t,r} \). The overall objective is to minimize the total system cost, composed of investment and operational components, as formalized in the following objective function:

\begin{flalign}\label{OB_Det_Problem}
&\min \enspace z = \displaystyle\sum_{i\in \Omega_S}  I^{pv}\gamma^{pv}_i + \displaystyle\sum_{i\in \Omega_S} \displaystyle\sum_{f\in \mathcal {F}} 
 I_{f}^{bt}\gamma^{bt}_{i,f} + 
\displaystyle\sum_{i\in \Omega_S} \displaystyle\sum_{t\in {\mathcal {T}}}  \displaystyle\sum_{r\in {\mathcal {R}}}  (\lambda^{bg}_{t,r} p\kappa^{bg}_{i,t,r} - \lambda^{sg}_{t,r} p\kappa^{sg}_{i,t,r})
\end{flalign}

\subsection{Network constraints}
The operation of the EC must comply with the physical constraints of the DN, including nodal power balance, voltage levels, and thermal limits on distribution lines. To preserve tractability in the optimization model, the network is represented using a second-order cone relaxation of the branch flow equations. The resulting constraints are detailed below.
\begin{subequations}
\begin{align}
&\displaystyle\sum_{(i,j) \in {\mathcal{L}}} p^{i,j}_{t,r} - \displaystyle\sum_{(j,i) \in {\mathcal{L}}} (p^{j,i}_{t,r} - R_{j,i}\ell^{j,i}_{t,r}) =    \left \{
  \begin{aligned}
    & \Delta p_{i,t,r}  && \text{if } i \in \Omega_A\\
    & \Delta \kappa p_{i,t,r} && \text{if } i \in \Omega_S\\
    & 0 && Otherwise
  \end{aligned} \right.\forall i \in \Omega ,\forall t \in {\mathcal{T}}, \forall r \in {\mathcal{R}}\label{Act_Node_Eq}
\end{align}
\begin{align}
&\displaystyle\sum_{(i,j) \in {\mathcal{L}}} q^{i,j}_{t,r} - \displaystyle\sum_{(j,i) \in {\mathcal{L}}} (q^{j,i}_{t,r} - X_{j,i}\ell^{j,i}_{t,r}) =   qg_{i,t,r} - \frac{1}{SB}QL_{i,t,r} && \forall i \in \Omega, \forall t \in {\mathcal{T}}, \forall r \in {\mathcal{R}} \label{React_Node_Eq}\\
&v_{j,t,r} = v_{i,t,r} -2 (R_{i,j} p^{i,j}_{t,r} + X_{i,j}q^{i,j}_{t,r}) + (R^2_{i,j} + X^2_{i,j})\ell^{i,j}_{t,r} && \forall t \in {\mathcal{T}}, \forall r \in {\mathcal{R}}, \forall (i,j) \in {\mathcal{L}} \label{Act_Line_Eq}  \\
&(p^{i,j}_{t,r})^2 + (q^{i,j}_{t,r})^2 \leq \ell^{i,j}_{t,r} v_{i,t,r} && \forall t \in {\mathcal{T}}, \forall r \in {\mathcal{R}}, \forall (i,j) \in {\mathcal{L}}\label{React_Line_Eq} \\
&(V^{min}_i)^2 \leq v_{i,t,r} \leq (V^{max}_i)^2 && \forall i \in \Omega, \forall t \in {\mathcal{T}}, \forall r \in {\mathcal{R}}\label{Eq_v_limit} \\
&p\kappa^{bg}_{i,t,r} \leq \frac{1}{SB}PBG^{max}_i && \forall i \in \Omega_S, \forall t \in {\mathcal{T}}, \forall r \in {\mathcal{R}} \label{KappaBG_C}\\
&p\kappa^{sg}_{i,t,r} \leq \frac{1}{SB}PSG^{max}_i && \forall i \in \Omega_S, \forall t \in {\mathcal{T}}, \forall r \in {\mathcal{R}} \label{KappaSG_C}\\
&-\frac{1}{SB}PSG^{max}_i \leq qg_{i,t,r} \leq \frac{1}{SB}PSG^{max}_i && \forall i \in \Omega, \forall t \in {\mathcal{T}}, \forall r \in {\mathcal{R}}\label{PSG_Reac_C}\\
&\Delta p_{i,t,r} = pv_{i,t,r} - \frac{1}{SB}PL_{i,t,r} + ds_{i,t,r} - ch_{i,t,r} &&\forall i \in \Omega_A, \forall t \in {\mathcal{T}}, \forall r \in {\mathcal{R}} \label{Eq_DeltaP}\\
&\Delta  \kappa p_{i,t,r} = p\kappa^{bg}_{i,t,r} - p\kappa^{sg}_{i,t,r} + pv^{cm}_{i,t,r} + \sum_{f \in \mathcal{F}} \left( ds^{cm}_{i,t,r,f} - ch^{cm}_{i,t,r,f} \right) && \forall i \in \Omega_S ,\forall t \in {\mathcal{T}}, \forall r \in {\mathcal{R}}\label{Eq_DeltaKappaP}
\end{align}
\end{subequations}

Equations~\eqref{Act_Node_Eq} and~\eqref{React_Node_Eq} enforce the active and reactive power balance at each node \( i \in \Omega \) for every time step \( t \in \mathcal{T} \) and year \( r \in \mathcal{Y} \). The left-hand side represents the net incoming active or reactive power, including losses through the terms \( R_{j,i} \ell^{j,i}_{t,r} \) and \( X_{j,i} \ell^{j,i}_{t,r} \), respectively, while \( \ell^{i,j}_{t,r} \) represents the squared magnitude of the current flowing through line \( (i,j) \). The right-hand side captures the net injection or withdrawal at each node. Thus, if the node belongs to an active user (\( i \in \Omega_A \)), it corresponds to the local net injection \( \Delta p_{i,t,r} \); if the node is the point of common coupling with the upstream grid (\( i \in \Omega_S \)), it corresponds to \( \Delta \kappa p_{i,t,r} \); and otherwise, the net injection is zero. For reactive power, it is assumed that each node may inject or absorb an amount \( qg_{i,t,r} \) in response to the reactive power requirement \( QL_{i,t,r} \). Equation~\eqref{Act_Line_Eq} models the voltage drop along each distribution line \( (i,j) \in \mathcal{L} \), accounting for the active and reactive power flows \( p^{i,j}_{t,r} \) and \( q^{i,j}_{t,r} \), the line resistance \( R_{i,j} \), and reactance \( X_{i,j} \). The voltage variable \( v_{i,t,r} \) denotes the squared voltage magnitude at the sending bus \( i \), and the expression reflects the impact of power flows and current-induced losses on the receiving-end voltage \( v_{j,t,r} \).

  The second-order cone constraint~\eqref{React_Line_Eq} ensures that the apparent power flowing through each line remains within the thermal limit, by enforcing that the squared magnitude of power flow does not exceed the product of current magnitude and sending-end voltage, where voltage magnitudes are bounded in Equation~\eqref{Eq_v_limit}. The maximum amount of active power that can be exchanged with the grid is limited by the rated capacity of the connection point in ~\eqref{KappaBG_C}--\eqref{KappaSG_C}, while ~\eqref{PSG_Reac_C} bounds the reactive power. Equation~\eqref{Eq_DeltaP} defines the net power injection \( \Delta p_{i,t,r} \) for each active user node \( i \in \Omega_A \) as the difference between PV generation \( pv_{i,t,r} \), scaled load \( PL_{i,t,r} \), and the charge/discharge operations of local storage systems. Specifically, \( ds_{i,t,r} \) and \( ch_{i,t,r} \) represent the discharging and charging power at node \( i \), respectively, during time \( t \) and year \( r \). Equation~\eqref{Eq_DeltaKappaP} characterizes the net injection at the point of common coupling with the grid. This term \( \Delta \kappa p_{i,t,r} \) includes the balance between power purchased and sold to the grid (\( p\kappa^{bg}_{i,t,r} \), \( p\kappa^{sg}_{i,t,r} \)), the generation of the shared PV system \( pv^{cm}_{i,t,r} \), and the aggregated charging/discharging operations of the community BESS. The latter is captured by the difference between \( ds^{cm}_{i,t,r,f} \) and \( ch^{cm}_{i,t,r,f} \) for each technology \( f \in \mathcal{F} \).

\subsection{Energy trading constraints}
This set of constraints defines the internal market structure of the EC and governs the logic of P2P energy exchanges, interactions with the external grid, and the utilization of community assets. 

\begin{subequations} \label{SetEq_Delta}
\begin{align}
	&\Delta p_{i,t,r} = \Delta p^{+}_{i,t,r} - \Delta p^{-}_{i,t,r} &&\forall i \in \Omega_A, \forall t \in {\mathcal{T}}, \forall r \in {\mathcal{R}} \label{Eq_DeltaTotal}\\
	&\Delta p^{+}_{i,t,r} = p^{sm}_{i,t,r} + p^{sg}_{i,t,r} + ch^{cm'}_{i,t,r} &&\forall i \in \Omega_A, \forall t \in {\mathcal{T}}, \forall r \in {\mathcal{R}} \label{Eq_DeltaMas} \\
	&\Delta p^{-}_{i,t,r} = p^{bm}_{i,t,r} + p^{bg}_{i,t,r} + ds^{cm'}_{i,t,r} + pv^{cm'}_{i,t,r} &&\forall i \in \Omega_A, \forall t \in {\mathcal{T}}, \forall r \in {\mathcal{R}} \label{Eq_DeltaMenos}\\
	&\Delta p^{+}_{i,t,r} \leq M (y_{i,t,r}) &&\forall i \in \Omega_A, \forall t \in {\mathcal{T}}, \forall r \in {\mathcal{R}} \label{Eq_y1} \\
	&\Delta p^{-}_{i,t,r} \leq M(1-y_{i,t,r}) &&\forall i \in \Omega_A, \forall t \in {\mathcal{T}}, \forall r \in {\mathcal{R}} \label{Eq_y2} \\
    &\displaystyle\sum_{i \in \Omega_A} p^{sm}_{i,t,r} = \displaystyle\sum_{i \in \Omega_A} p^{bm}_{i,t,r} &&\forall t \in {\mathcal{T}}, \forall r \in {\mathcal{R}} \label{Eq_Balance_1} \\
	&\displaystyle\sum_{i \in \Omega_A} p^{sg}_{i,t,r} + \displaystyle\sum_{i \in \Omega_S} ds^{cm'}_{i,t,r} + \displaystyle\sum_{i \in \Omega_S} pv^{cm'}_{i,t,r}  \geq \displaystyle\sum_{i \in \Omega_S} p\kappa^{sg}_{i,t,r} &&\forall t \in {\mathcal{T}}, \forall r \in {\mathcal{R}} \label{Eq_sold_grid} \\
	&\displaystyle\sum_{i \in \Omega_A} p^{bg}_{i,t,r} + \displaystyle\sum_{i \in \Omega_S} ch^{cm'}_{i,t,r} \leq \displaystyle\sum_{i \in \Omega_S} p\kappa^{bg}_{i,t,r} &&\forall t \in {\mathcal{T}}, \forall r \in {\mathcal{R}} \label{Eq_bought_grid}
\end{align}
\end{subequations}
Equation~\eqref{Eq_DeltaTotal} decomposes the net injection \( \Delta p_{i,t,r} \) of each user \( i \in \Omega_A \) into a positive component \( \Delta p^{+}_{i,t,r} \) (energy surplus) and a negative component \( \Delta p^{-}_{i,t,r} \) (energy deficit). These components are further detailed in Equations~\eqref{Eq_DeltaMas}--\eqref{Eq_DeltaMenos}, where surplus may be it can be sold to a neighbor via the internal market (\( p^{sm}_{i,t,r} \)), exported to the grid (\( p^{sg}_{i,t,r} \)), or allocated to the community BESS (\( ch^{cm'}_{i,t,r} \)). Conversely, a deficit can be met by buying from a neighbor (\( p^{bm}_{i,t,r} \)), purchasing from the grid (\( p^{bg}_{i,t,r} \)), discharging energy from the community BESS (\( ds^{cm'}_{i,t,r} \)), or using part of the shared PV generation (\( pv^{cm'}_{i,t,r} \)). To prevent simultaneous buying and selling by the same user, a binary variable \( y_{i,t,r} \) is introduced. Equations~\eqref{Eq_y1}--\eqref{Eq_y2} enforce mutual exclusivity by activating only one side of the transaction per time step, where \( M \) is a sufficiently large number.

  The internal market is cleared locally at each time \( t \) and year \( r \) by Equation~\eqref{Eq_Balance_1}, which guarantees that the total amount of energy sold among peers (\( p^{sm}_{i,t,r} \)) matches the total amount purchased (\( p^{bm}_{i,t,r} \)). External grid interactions are bounded by Equations~\eqref{Eq_sold_grid}--\eqref{Eq_bought_grid}. The total amount of energy sold to the grid must not exceed the aggregated surplus available at the PCC, including contributions from individual sales, community battery discharges, and shared PV generation. Similarly, the energy purchased from the grid must at least cover the aggregate demand that cannot be satisfied internally or via the community battery. These constraints are expressed as inequalities to account for minor losses in the distribution lines, which prevent perfect balance between imported/exported energy and internal transactions.

\subsection{User Assets}
This set of constraints models the operation of individually owned assets, including BESS and PV systems.
\begin{subequations}\label{Set_BT_User}
\begin{align}
\label{Eq_SOC}  	
	&soc_{i,t,r}= soc_{i,t-1,r}+\varphi^{ch} ch_{i,t,r}-\frac{1}{\varphi^{ds} }ds_{i,t,r} &&\forall i \in \Omega_A, \forall t \in {\mathcal{T}}, \forall r \in {\mathcal{R}}\\
\label{Eq_SOC_limit} 
	&SOC^{min} \frac{1}{SB}\Gamma^{bt}_i \leq soc_{i,t,r} \leq SOC^{max} \frac{1}{SB}\Gamma^{bt}_i &&\forall i \in \Omega_A, \forall t \in {\mathcal{T}}, \forall r \in {\mathcal{R}}\\
\label{Eq_ch}
	&ch_{i,t,r}  \leq \frac{1}{SB}PB (w_{i,t,r}) &&\forall i \in \Omega_A, \forall t \in {\mathcal{T}}, \forall r \in {\mathcal{R}}\\
\label{Eq_ds} 
	&ds_{i,t,r} \leq \frac{1}{SB}PB (1-w_{i,t,r}) - \frac{1}{SB}PB(1-BT_i) &&\forall i \in \Omega_A, \forall t \in {\mathcal{T}}, \forall r \in {\mathcal{R}}\\
\label{bin_bat_cap}
    &w_{i,t,r} \leq BT_i &&\forall i \in \Omega_A, \forall t \in {\mathcal{T}}, \forall r \in {\mathcal{R}}\\
&pv_{i,t,r}\leq PG^{max}_{t,r} \frac{1}{SB}\Gamma^{pv}_i && \forall i \in \Omega_A, \forall t \in {\mathcal{T}}, \forall r \in {\mathcal{R}} \label{Eq_pg_limit}
\end{align}
\end{subequations}  
Equation~\eqref{Eq_SOC} defines the state-of-charge dynamics of the user's battery, considering charging \( ch_{i,t,r} \) and discharging \( ds_{i,t,r} \) operations, scaled by their respective efficiencies \( \varphi^{ch} \) and \( \varphi^{ds} \). The stored energy is bounded by Equation~\eqref{Eq_SOC_limit}, according to the installed capacity \( \Gamma^{bt}_i \) and the state-of-charge limits \( SOC^{min} \) and \( SOC^{max} \).  Charging and discharging power are limited by \( PB \), the battery’s maximum rate, and controlled via the binary variable \( w_{i,t,r} \), which prevents simultaneous actions, as implemented in Equations~\eqref{Eq_ch} and~\eqref{Eq_ds}. Equation~\eqref{bin_bat_cap} deactivates the battery if it is not installed (\( BT_i = 0 \)). Finally, Equation~\eqref{Eq_pg_limit} bounds PV generation based on the normalized availability profile \( PG^{max}_{t,r} \), scaled by the installed capacity \( \Gamma^{pv}_i \).

\subsection{Community Assets}
The operation of community PV and BESS assets follows the same logic as for individual devices:
\begin{subequations}\label{Set_BT_CB}
\small
\begin{align}
\label{Eq_SOC_CB}  	
	&soc^{cm}_{i,t,r,f}= soc^{cm}_{i,t-1,r,f}+\varphi^{ch} ch^{cm}_{i,t,r,f}-\frac{1}{\varphi^{ds} }ds^{cm}_{i,t,r,f} &&\forall i \in \Omega_S, \forall t \in {\mathcal{T}}, \forall r \in {\mathcal{R}}, \forall f \in {\mathcal{F}}\\
\label{Eq_SOC_limit_CB} 
	&SOC^{min} \gamma^{bt}_{i,f} DBT_{r,f} \leq soc^{cm}_{i,t,r,f} \leq SOC^{max} \gamma^{bt}_{i,f} DBT_{r,f} &&\forall i \in \Omega_S, \forall t \in {\mathcal{T}}, \forall r \in {\mathcal{R}}, \forall f \in {\mathcal{F}}\\
\label{Eq_ch_CB} 
	&ch^{cm}_{i,t,r,f}  \leq \frac{1}{SB}PB^{cm}(w^{cm}_{i,t,r,f}) &&\forall i \in \Omega_S, \forall t \in {\mathcal{T}}, \forall r \in {\mathcal{R}}, \forall f \in {\mathcal{F}}\\
\label{Eq_ds_CB} 
	&ds^{cm}_{i,t,r,f} \leq \frac{1}{SB}PB^{cm}(1-w^{cm}_{i,t,r,f}) -  \frac{1}{SB}PB^{cm}(1-\nu^{cm}_{f}) &&\forall i \in \Omega_S, \forall t \in {\mathcal{T}}, \forall r \in {\mathcal{R}}, \forall f \in {\mathcal{F}}\\
    &w^{cm}_{i,t,r,f} \leq \nu^{cm}_{f}  &&\forall i \in \Omega_S, \forall r \in {\mathcal{R}}, \forall f \in {\mathcal{F}} \label{Nu_Bess}\\
	&\gamma^{bt}_{i,f} \leq \nu^{cm}_{f} \frac{1}{SB}\Gamma^{max}  &&\forall i \in \Omega_S, \forall r \in {\mathcal{R}}, \forall f \in {\mathcal{F}} \label{Gamma_BESS_CM} \\
    &\displaystyle\sum_{f \in \mathcal{F}} \nu^{cm}_{f} \leq 1 && \label{BESS_Tech} \\
&pv^{cm}_{i,t,r}\leq PG^{max}_{t,r} \gamma^{pv}_i && \forall i \in \Omega_S, \forall t \in {\mathcal{T}}, \forall r \in {\mathcal{R}} \label{Eq_pg_Community}\\
&\gamma^{pv}_i  \leq \frac{1}{SB}\Gamma^{max}  &&\forall i \in \Omega_S \label{Max_PV_CM_Inst}\\
\label{Eq_ds_dscb} 
	&\displaystyle\sum_{i \in \Omega_B} ds^{cm'}_{i,t,r} = \displaystyle\sum_{i \in \Omega_S} \displaystyle\sum_{f \in \mathcal{F}}ds^{cm}_{i,t,r,f} &&\forall t \in {\mathcal{T}}, \forall r \in {\mathcal{R}}\\
\label{Eq_ch_chcb} 
	&\displaystyle\sum_{i \in \Omega_B} ch^{cm'}_{i,t,r} = \displaystyle\sum_{i \in \Omega_S} \displaystyle\sum_{f \in \mathcal{F}} ch^{cm}_{i,t,r,f} &&\forall t \in {\mathcal{T}}, \forall r \in {\mathcal{R}}\\
&\displaystyle\sum_{i \in \Omega_B} pv^{cm'}_{i,t,r} = \displaystyle\sum_{i \in \Omega_S} pv^{cm}_{i,t,r} &&\forall t \in {\mathcal{T}}, \forall r \in {\mathcal{R}}\label{Community_PV}
\end{align}
\end{subequations}
The subscript “cm” denotes variables associated with community-level assets, and the main differences lie in the modeling of degradation and technology selection. In particular, the state-of-charge bounds in Equation~\eqref{Eq_SOC_limit_CB} incorporate a degradation factor \( DBT_{r,f} \), which reduces the usable capacity \( \gamma^{bt}_{i,f} \) over time. The community can install at most one type of BESS technology, enforced by the binary variable \( \nu^{cm}_f \) in Equation~\eqref{BESS_Tech}. Equation~\eqref{Gamma_BESS_CM} further bounds the installed capacity per user and technology, while Equation~\eqref{Nu_Bess} restricts operation to selected technologies only. Similar capacity limits apply to PV systems, as defined in Equations~\eqref{Eq_pg_Community} and~\eqref{Max_PV_CM_Inst}. Finally, Equations~\eqref{Eq_ds_dscb}–\eqref{Community_PV} ensure consistency between the total energy exchanged with the community devices and the sum of individual contributions from users, linking centralized operation with user-level accounting.

\newpage
\section{Case study and computational results}\label{sec:results}

This section presents a comprehensive sensitivity analysis aimed at systematically evaluating the impact of variations in key parameters on system performance and battery installation decisions. The analysis begins with a description of the case study, followed by a detailed examination of the different scenarios considered.

\subsection{Case Study} 
\begin{figure}[H]
    \centering
    \includegraphics[width=0.5\textwidth]{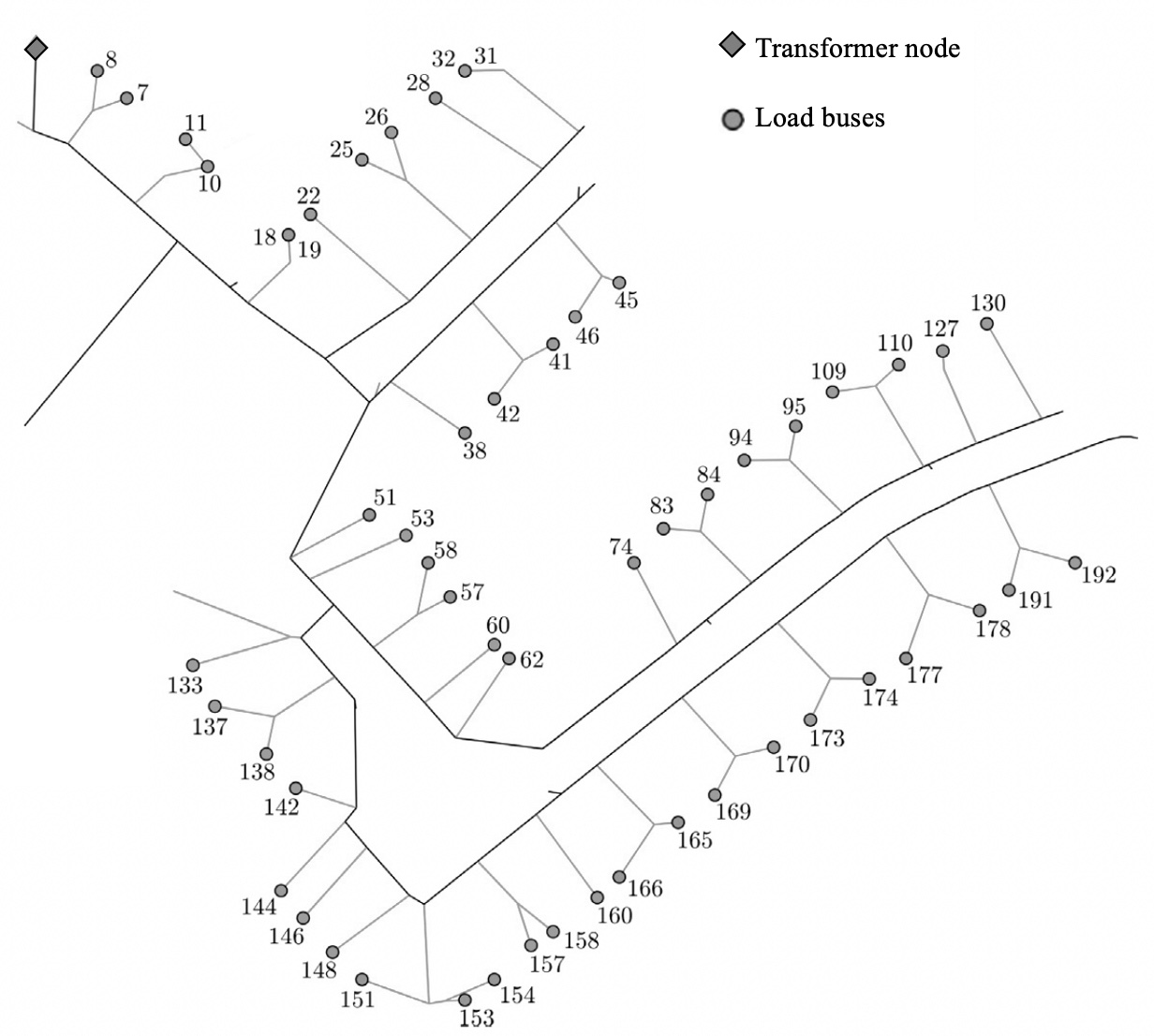} 
    \caption{Simplified 206-node low-voltage test feeder \cite{Garcia2024}}
    \label{fig:LV_DN_Case}
\end{figure}

The parameters used in this case study are based on the reduced model of the IEEE European Low Voltage Test Feeder proposed in \cite{Khan2022}, which achieves identical results with up to 80\% less computation time by simplifying the original network. Building on this approach, the present study adopts a further simplification strategy similar to that in \cite{Garcia2024}, where additional non-essential junction nodes are removed. As a result, the DN is reduced from 906 to 206 nodes as shown in Figure \ref{fig:LV_DN_Case}, further improving computational efficiency.

Table~\ref{tab:lv_network_summary} summarizes the key parameters of the case study. The network comprises 56 buses and 205 branches, serving 54 consumption points. Among these, 31 are prosumers equipped with PV systems, and 16 of them also have BESS installed, while the remaining 23 correspond to traditional consumers without any generation or storage assets. The energy community exhibits a total peak active demand of 0.374~MW and a peak PV generation of 0.195~MW, resulting in a PV penetration level of approximately 52\%. The total installed BESS capacity is 0.088~MWh, which corresponds to roughly 15\% of the community’s peak hourly demand.

\begin{table}[H]
\centering
\footnotesize
\caption{Key parameters of the case study}
\begin{tabular}{|l|c|c|}
\hline
\textbf{Parameter} & \textbf{Value} & \textbf{Unit} \\
\hline
\multicolumn{3}{|c|}{\textbf{Energy Community (EC)}} \\
\hline
Peak active power load (EC) & 0.374 & MW \\
Peak reactive power load (EC) & 0.0748 & MVar \\
Peak PV generation (EC) & 0.195 & MW \\
Total installed BESS capacity (EC) & 0.088 & MWh \\
\hline
\multicolumn{3}{|c|}{\textbf{Network Topology}} \\
\hline
Number of branches & 205 & -- \\
Number of buses & 56 & -- \\
Number of consumption points & 54 & -- \\
\hline
\multicolumn{3}{|c|}{\textbf{Per-User Values}} \\
\hline
Peak active power load  & 0.0068 & MW \\
Peak reactive power load  & 0.00136 & MVar \\
Peak PV generation  & 0.008 & MW \\
Peak BESS capacity  & 0.0055 & MWh \\
\hline
\end{tabular}\label{tab:lv_network_summary}
\end{table}

\begin{figure}[H]
    \centering
    \begin{minipage}{0.43\textwidth}
        \centering
        \includegraphics[width=\textwidth]{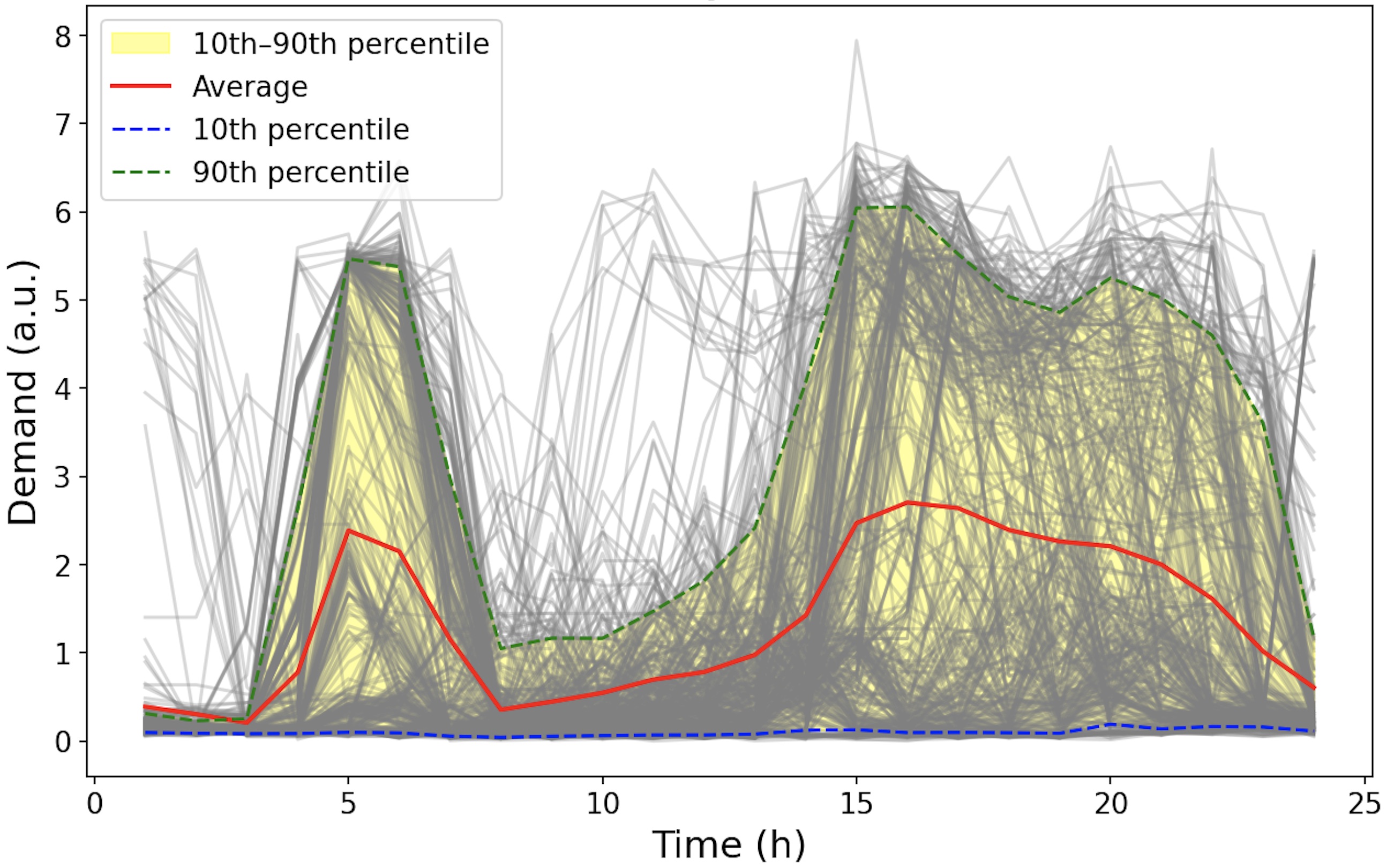}
        \caption{Daily demand profiles}
        \label{fig:demand_profiles}
    \end{minipage}
    \hfill
    \begin{minipage}{0.50\textwidth}
        \centering
        \includegraphics[width=\textwidth]{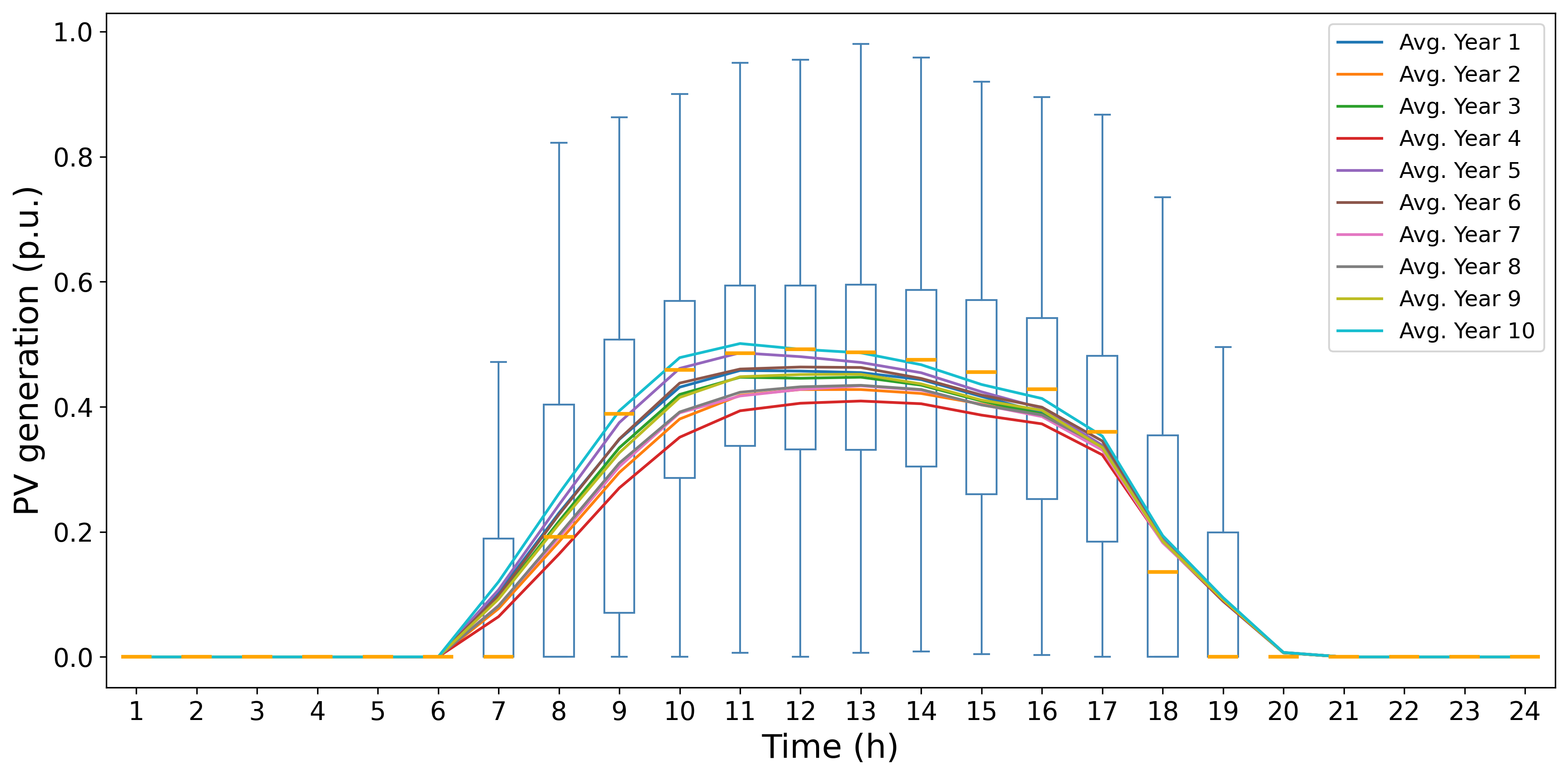}
        \caption{Daily generation profile}
        \label{fig:gen_prices}
    \end{minipage}
\end{figure}

This study considers representative rural households located in Spain. Figures~\ref{fig:demand_profiles} and~\ref{fig:gen_prices} present the daily electricity demand and photovoltaic generation profiles used in the analysis, based on typical values for this region and household type. Specifically, Figure~\ref{fig:demand_profiles} shows the daily electricity demand profiles. The model assumes an expected demand profile for each year (red line), which is computed as the average of 365 individual daily profiles generated for that year. These annual profiles are projected based on the historical data and methodology described in Section~\ref{sec:Methodology}, and then replicated to obtain expected demand profiles for the remaining nine years in the planning horizon. The shaded region represents the variability across days, while the blue and black dashed lines correspond to the 10th and 90th percentiles, respectively. Similarly, Figure~\ref{fig:gen_prices} presents the expected PV generation profiles for each year (colored lines), together with boxplots that capture the hourly variability across a 10-year horizon. Each boxplot aggregates 3650 values per hour (365 days $\times$ 10 years), showing the distribution of PV generation over time. This representation captures both the temporal and inter-annual variability of solar production. The colored curves reflect the expected generation profile for each year, which is used as the deterministic input for the optimization model.

\begin{figure}[H]
    \centering
    \begin{minipage}{0.45\textwidth}
        \centering
        \includegraphics[width=\textwidth]{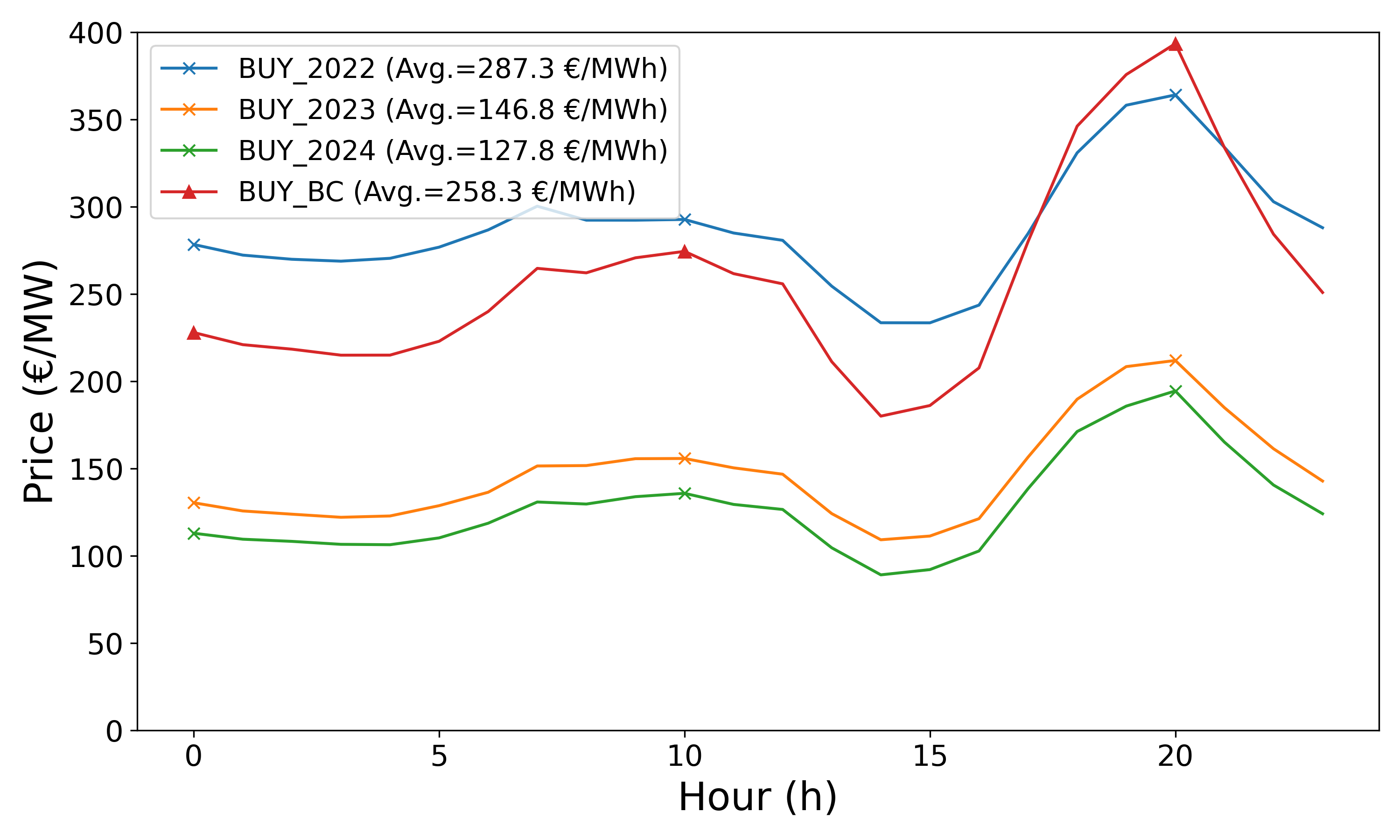}
        \caption{Hourly average buying prices for the last three years}
        \label{fig:buying_prices}
    \end{minipage}
    \hspace{0.05\textwidth}
    \begin{minipage}{0.45\textwidth}
        \centering
        \includegraphics[width=\textwidth]{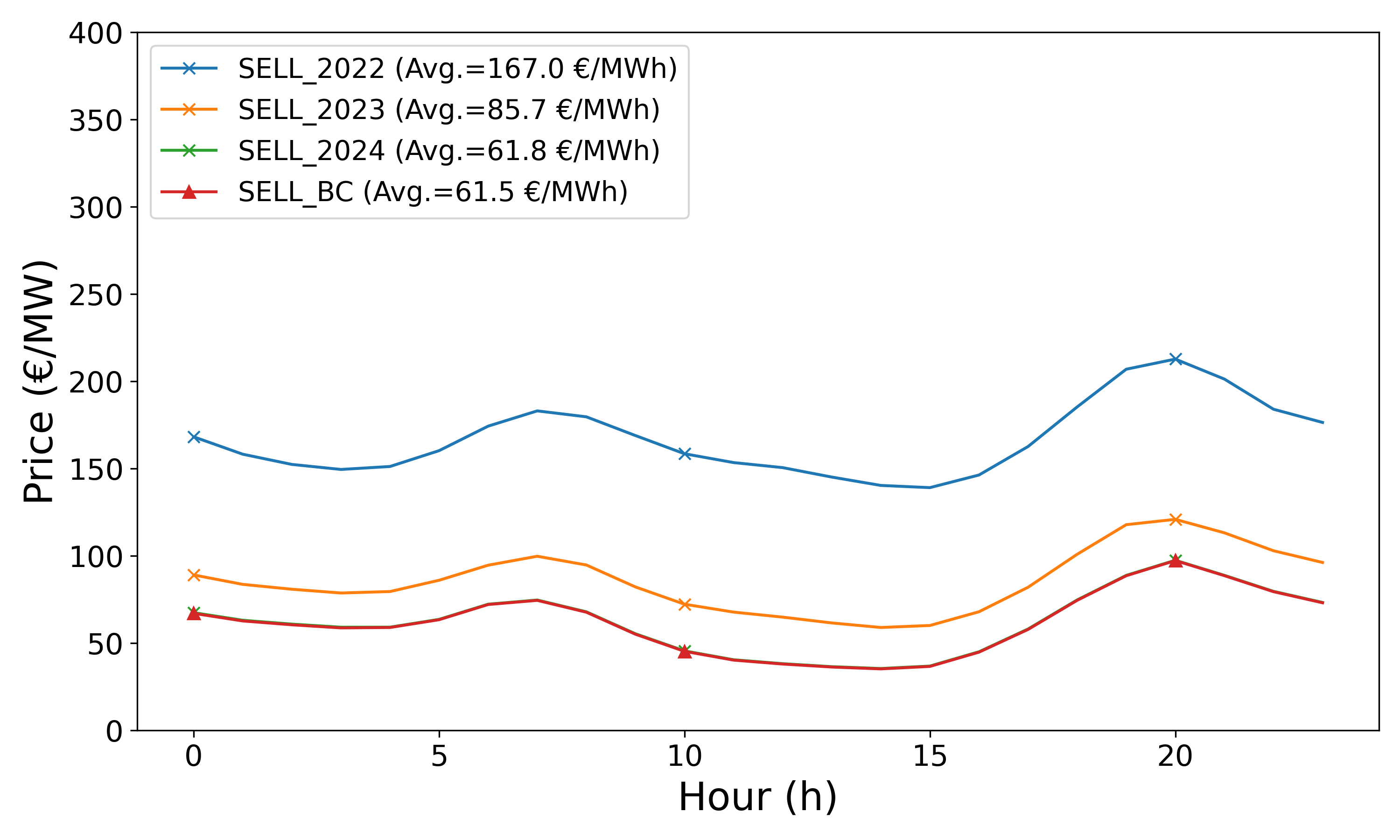}
        \caption{Hourly average selling prices for the last three years}
        \label{fig:selling_prices}
    \end{minipage}
\end{figure}

Regarding the profile for electricity grid prices, the electricity market is characterized by significant fluctuations, and this is consistent with the findings of~\cite{Wan2025}, which also highlight an anticipated upward trend in electricity prices. However, the data shown in Figures~\ref{fig:buying_prices} and~\ref{fig:selling_prices} reveal a decreasing trend in hourly market prices over the last three years (2022--2024). In particular, the average buying price dropped from \euro{}287.3/MWh in 2022 to \euro{}127.8/MWh in 2024, while the selling price fell from \euro{}167.0/MWh to \euro{}61.8/MWh over the same period. This reduction also led to a gap margin between buying and selling prices, thereby decreasing the potential arbitrage opportunities in the energy market. These average values are explicitly indicated in the legends of each curve, providing an overview of the annual variations. Given lower electricity prices, especially buying prices, the base case price scenario for this study adopts a more favorable economic context. Specifically, it assumes a buying price profile closer to the 2022 and adjusted with a 2\% annual increase to reflect projected trends. For the selling price, the 2024 profile is retained as a conservative estimate. This approach is justified by the need to maintain a sufficient spread between purchase and sale prices to incentivize the deployment of BESS, which, while increasingly competitive, are still sensitive to the price signals of the electricity market. A more detailed economic analysis of this impact is presented in a later section.

\begin{figure}[H]
    \centering
    \begin{minipage}{0.45\textwidth}
        \centering
        \includegraphics[width=\textwidth]{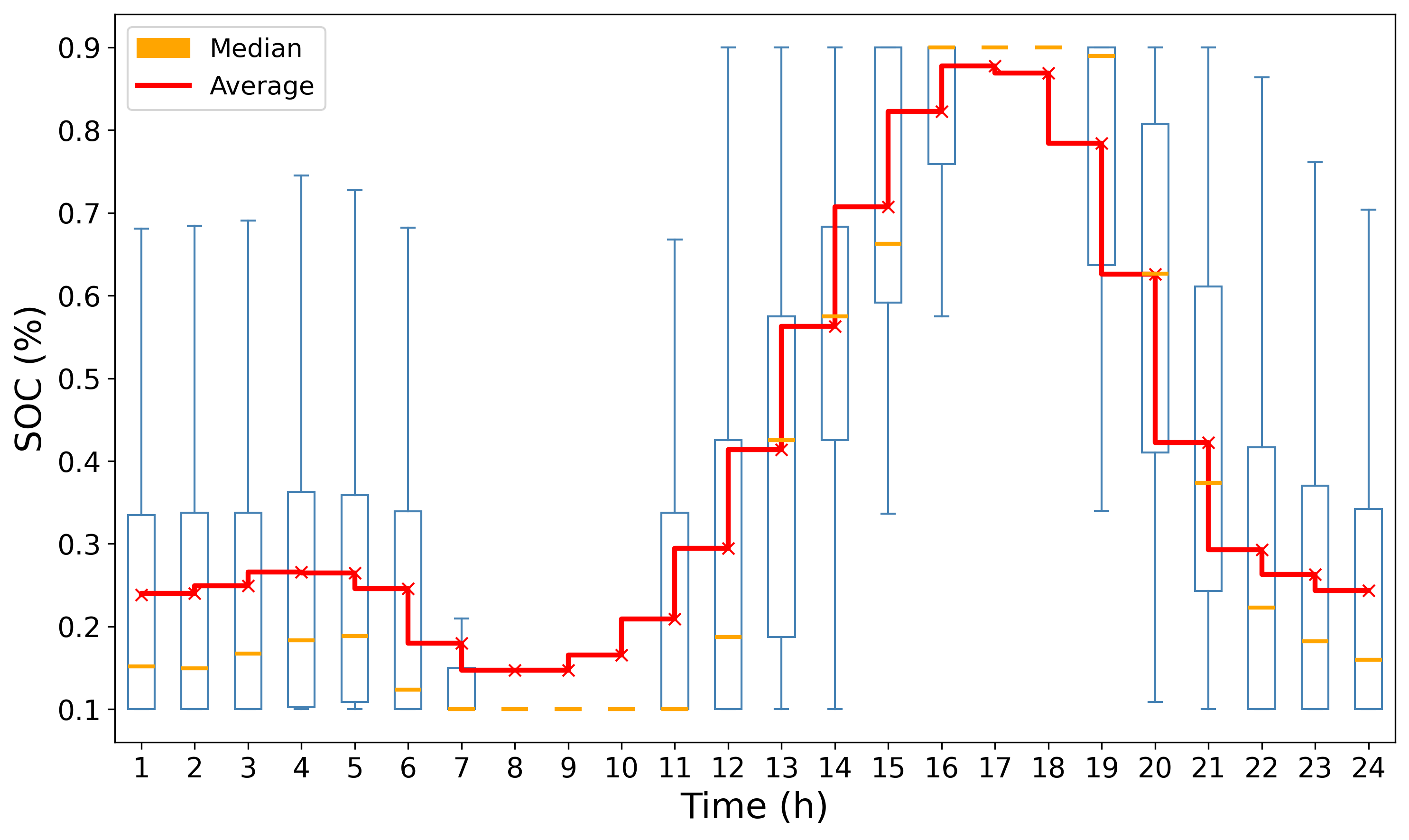}
        \caption{SOC of the battery}
        \label{fig:soc}
    \end{minipage}
    \hspace{0.05\textwidth}
    \begin{minipage}{0.45\textwidth}
        \centering
        \includegraphics[width=\textwidth]{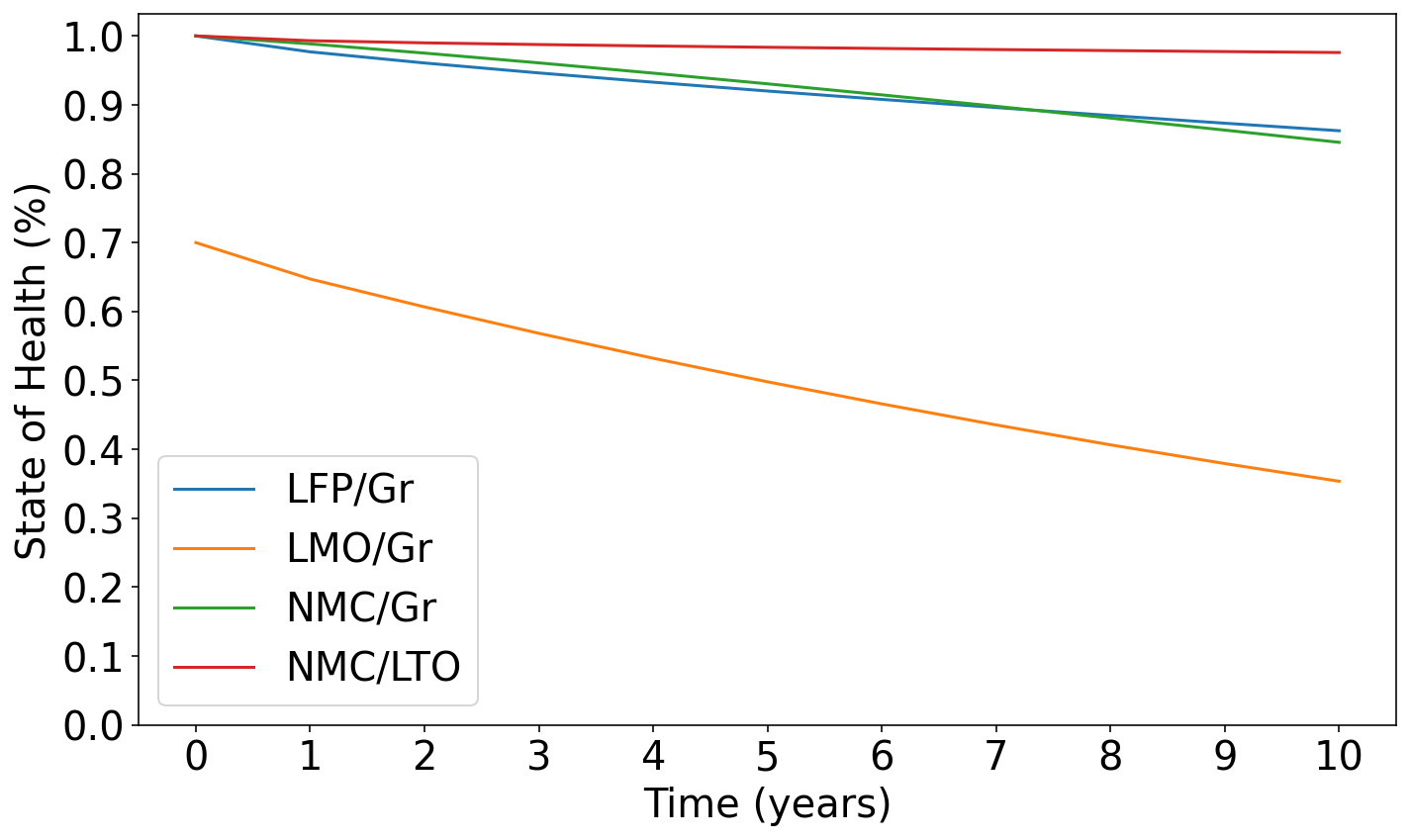}
        \caption{Degradation profiles}
        \label{fig:deg_bat}
    \end{minipage}
\end{figure}

Following the methodology described in Section~\ref{sec:Methodology}, the first step to estimate the depreciation of different BESS technologies consists of obtaining the SOC profile of a typical user. This profile is derived from the user's deterministic optimization model and is shown in Figure~\ref{fig:soc}.

As seen in the figure, the red line represents the average hourly SOC over 365 days, while the orange line shows the median value, and the boxplots depict the variability due to seasonal demand and generation patterns. A clear daily cycle emerges: the battery charges during daylight hours, corresponding to PV generation, and discharges during the evening when electricity prices are highest. Additionally, the battery also tends to charge slightly during the early morning hours when electricity prices are typically low. The average SOC remains higher than the median during the night and morning hours, indicating that on many days the battery is only partially charged, while on fewer days it reaches higher SOC levels, hence the skew. Despite the observed variability, the operational pattern is consistent. Two distinct behaviors stand out: (i) the battery tends to be fully discharged around 8--10~h, and (ii) it reaches nearly full capacity around 17--18~h. It is important to note that, since the optimization model does not enforce a fixed terminal SOC at the end of each 24-hour period, the profile naturally evolves from day to day. However, across the simulation horizon of one year (8760 hours), the battery tends to stabilize within a 20--30\% SOC range at the end of each day.

Using this SOC profile as input, a degradation simulation was performed in MATLAB for different BESS technologies, and the resulting evolution of the SoH over a 10-year horizon is shown in Figure~\ref{fig:deg_bat}. The results are consistent with expectations: the SL battery (LFP/Gr) exhibits the fastest degradation, with its SoH dropping significantly over the 10 years. Once the battery reaches a SOH of 40\%, it is regarded as depleted and consequently unsuitable for further use  \cite{Montoya}. The LMO/Gr and NMC/Gr technologies display similar degradation rates, remaining above 80\% SoH. In contrast, the NMC/LTO battery demonstrates minimal degradation, maintaining a nearly constant SoH close to 100\%, highlighting its superior cycling stability.


Regarding BESS price assumptions, the values adopted in this study are based on prevailing commercial prices from current market sources. These correspond to \textbf{Case~C1} in Table~\ref{table:battery_prices}, and represent the initial benchmark for all comparisons. In addition, three alternative price scenarios (\textbf{C2} to \textbf{C4}) are introduced to reflect expected cost reductions in the coming years, in line with projections reported in~\cite{Cheng2022} for the year 2030. Specificall, \textbf{C2} considers a moderate reduction of approximately 30\% from current prices. \textbf{C3} assumes a 40\% reduction, representing a more optimistic scenario, and \textbf{C4} corresponds directly to the projected 2030 values cited in~\cite{Cheng2022}. For consistency, all reductions are applied proportionally to each battery chemistry. The base chemistry selected for these projections is LFP/Gr, which is the most widely used in practice. Note that the LMO/Gr battery represents a SL battery technology, and its pricing reflects its lower cost and remaining capacity. In our analysis, SL batteries are assumed to begin operational at 70\% SoH, a value commonly cited as the minimum for their reuse \cite{Ma2024}. All first-life batteries (LFP/Gr, NMC/Gr, NMC/LTO) start at 100\% SoH.

\begin{table}[H]
\centering
\footnotesize
\caption{BESS price scenarios (\euro{}/kWh) used for the sensitivity analysis}
\begin{tabular}{|c|c|c|c|c|}
\hline
\textbf{Case} & \textbf{LFP/Gr (FL)} & \textbf{LMO/Gr (SL)} & \textbf{NMC/Gr (FL)} & \textbf{NMC/LTO (FL)} \\
\hline
\textbf{C1: Current commercial prices} & 685 & 277 & 987 & 3736 \\
\textbf{C2: Projected --30\%}          & 480 & 190 & 690 & 2610 \\
\textbf{C3: Projected --40\%}          & 410 & 160 & 590 & 2240 \\
\textbf{C4: Projected 2030 \cite{Cheng2022}} & 150 & 60  & 215 & 810 \\
\hline
\end{tabular}
\label{table:battery_prices}
\end{table}



\subsection{Global results}


Considering the structure of the studied energy community, where peer-to-peer energy trading is already in place and the penetration of PV systems and BESS reaches 52\% and 15\% of the community's peak electricity demand, respectively. The goal of the optimization model is to assess whether it is economically viable to complement these existing individual assets with a shared, community-scale PV-BESS installation that further enhances collective autonomy and operational flexibility. The first relevant insight from the base case analysis is that, under current or recent electricity market conditions (i.e., years 2023 and 2024), and assuming commercial BESS prices (scenario C1) for both FL and SL technologies, the installation of community PV systems is consistently justified. In contrast, the deployment of community BESS remains economically infeasible. This outcome is primarily due to the high capital cost of storage, which still represents a significant barrier, especially under moderate or declining electricity price levels. Notably, even SL batteries, despite their reduced investment cost, are not economically viable under current conditions. Consequently, BESS installation only becomes viable when either BESS prices are sufficiently low or electricity price spreads are high enough to enable profitable arbitrage. It is important to note that, in this study, the number of prosumers, the demand profile, and the network restrictions were assumed to remain constant.

\begin{figure}[H]
    \centering
    \begin{minipage}{0.45\textwidth}
        \centering
        \includegraphics[width=\textwidth]{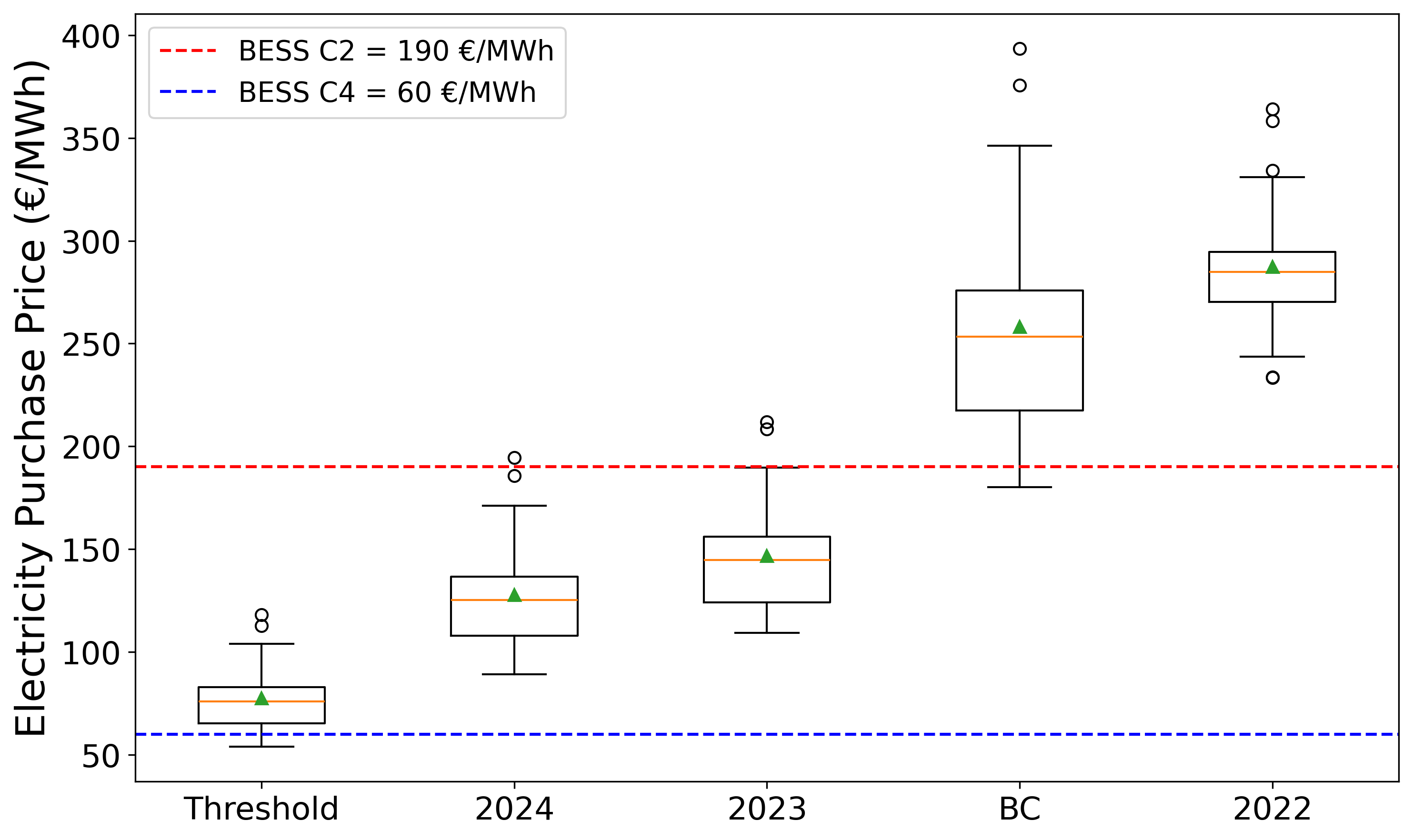}
        \caption{Electricity purchase price (2022–2024)}
        \label{fig:boxplot_1}
    \end{minipage}
    \hfill
    \begin{minipage}{0.45\textwidth}
        \centering
        \includegraphics[width=\textwidth]{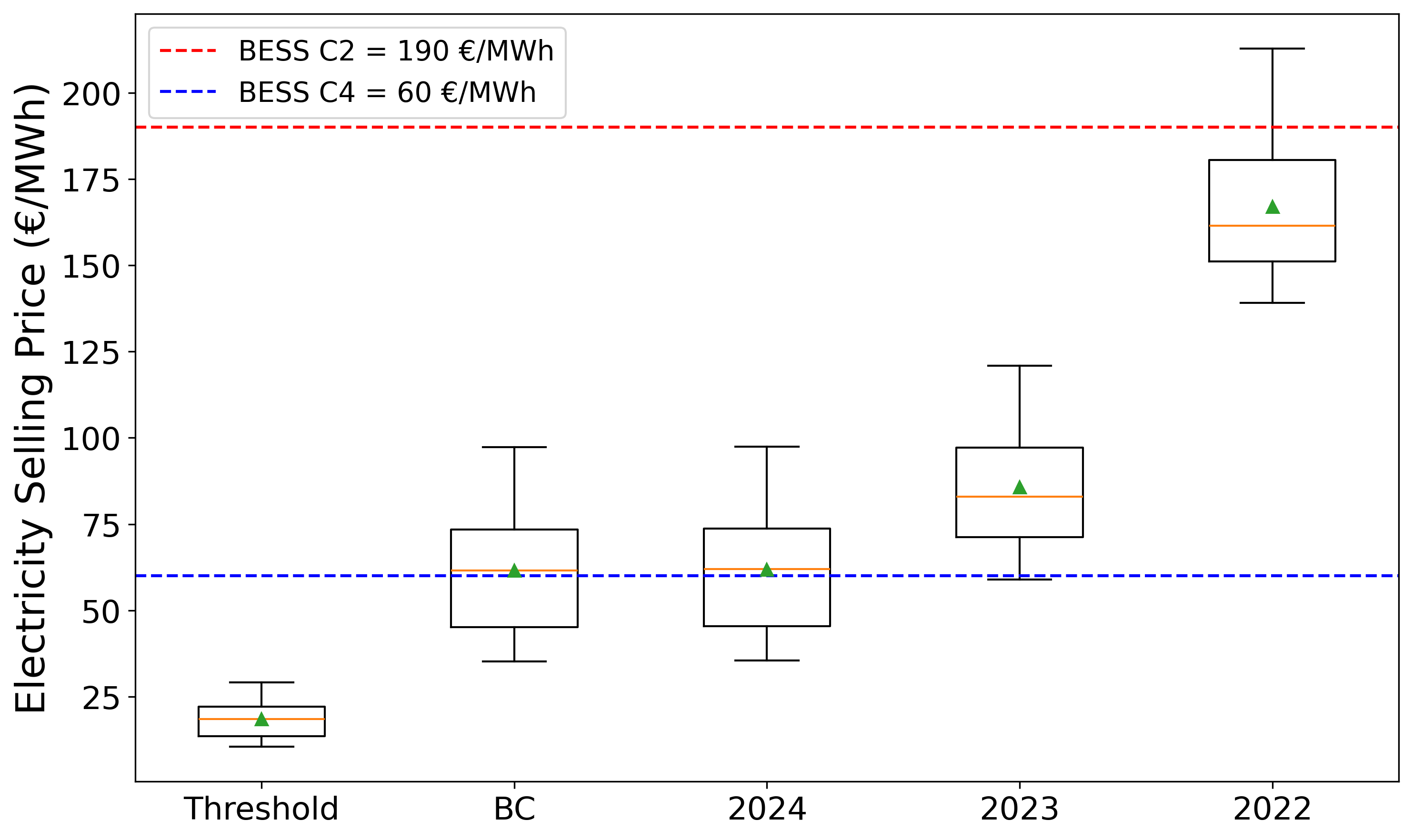}
        \caption{Electricity selling price (2022–2024)}
        \label{fig:boxplot_2}
    \end{minipage}
\end{figure}

To better understand the conditions under which BESS investment becomes feasible, a sensitivity analysis was conducted, focusing on the interaction between electricity market prices and SL battery costs. Figure~\ref{fig:boxplot_1} and Figure~\ref{fig:boxplot_2} present the boxplots of electricity purchase and selling prices in Spain for the years 2022 to 2024, as previously introduced in the case study. These plots capture the temporal variability of hourly electricity prices across the three-year period. Additionally, the dashed lines in each figure represent the cost of the SL BESS technologies under different scenarios. Specifically, the red dashed line indicates the SL BESS cost associated with scenario C2, while the blue line reflects the reduced investment level considered in scenario C4, both as specified in Table~\ref{table:battery_prices}. From these figures, it becomes evident that achieving economic viability for battery investments strongly depends on the combination of electricity prices and battery cost assumptions. For instance, under cost scenario C2 (\euro{}190/MWh for SL BESS), electricity purchase prices must at least reach the average values observed in 2022 for the installation of BESS to be profitable. In contrast, if SL battery costs decrease to scenario C4 (\euro{}60/MWh), then installation becomes viable even under the much lower electricity prices of 2024 corresponding to the threshold electricity purchase price case. Notably, in both cases where investment becomes favorable, the technology selected is an SL battery. This suggests that, at this stage, the investment cost is a more decisive factor than degradation in determining the feasibility of community-scale storage deployment.

\subsection{Battery price and degradation}\label{subsec:batprice/deg}

\begin{figure}[H]
    \centering
    \begin{minipage}{0.48\textwidth}
        \centering
        \includegraphics[width=\textwidth]{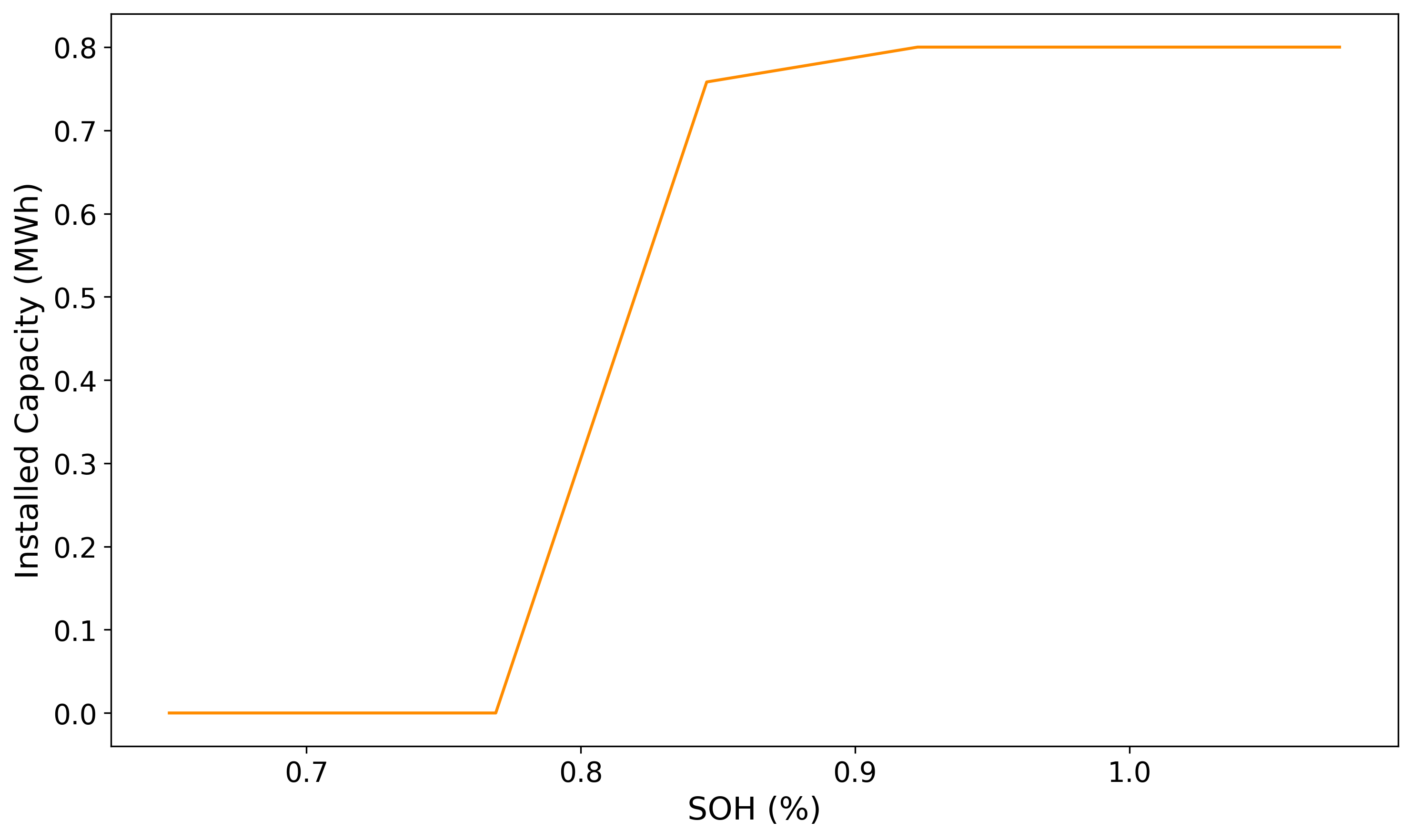}
        \caption{Installed SL BESS capacity as a function of State of Health (SOH)}
        \label{fig:deg}
    \end{minipage}
    \hfill
    \begin{minipage}{0.48\textwidth}
        \centering
        \includegraphics[width=\textwidth]{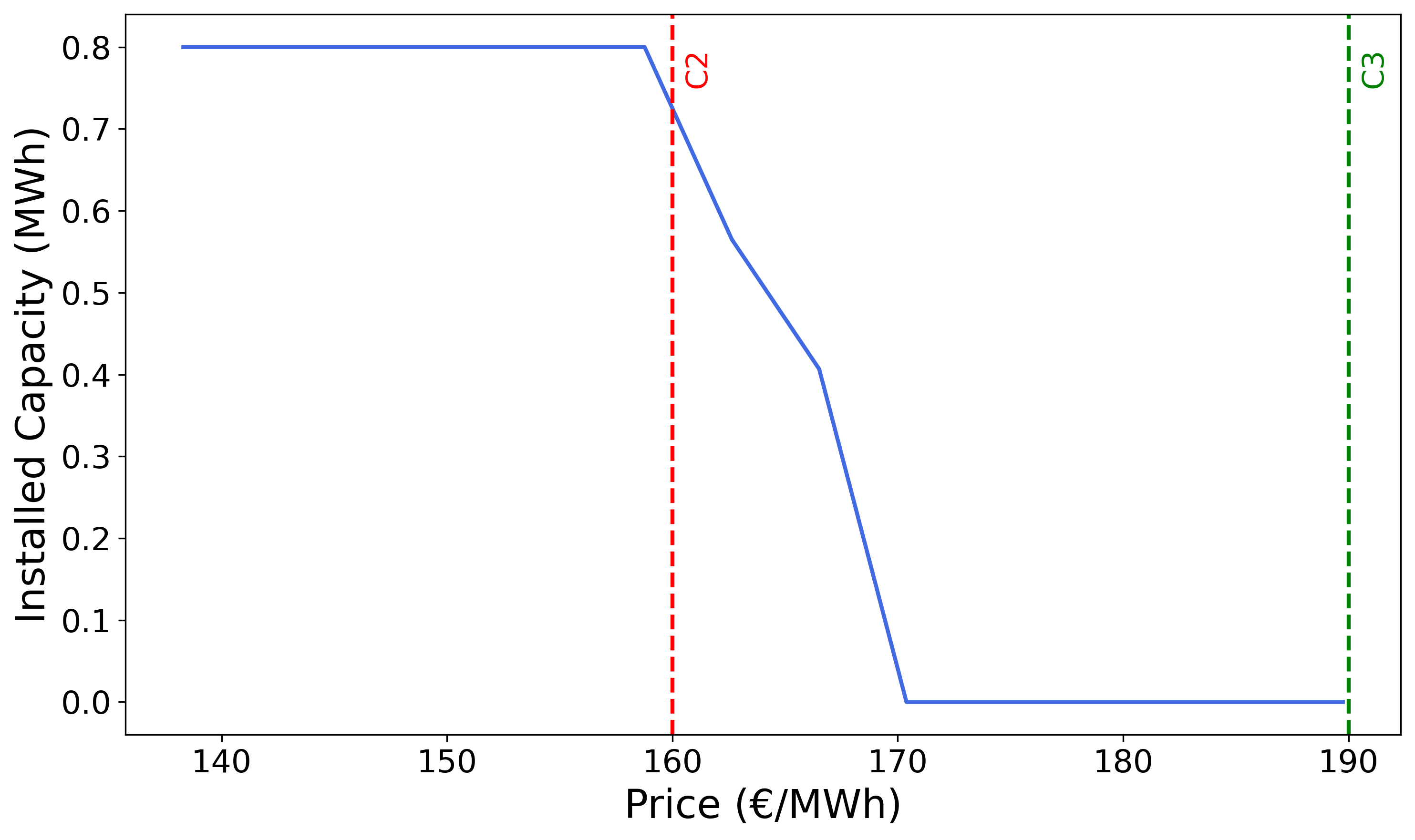}
        \caption{Installed SL BESS capacity as a function of SL BESS price}
        \label{fig:cp}
    \end{minipage}
\end{figure}

While the global results highlighted electricity prices and battery investment costs as the main drivers of BESS deployment, degradation also plays a critical role in the installation decision. To isolate its effect, a detailed analysis was conducted focusing on a representative scenario using SL BESS at a cost of \euro{}190/MWh (scenario C2), under the base case electricity price profile. In this configuration, when degradation is not considered, the model consistently selects to install the battery. However, once degradation is explicitly included using the current degradation profile of second-life technologies, the installation becomes economically unviable—despite maintaining the same cost and market conditions. This contrast reveals that degradation alone can revert the investment decision. To explore whether improvements in second-life battery performance could reverse this outcome, a series of simulations were performed by gradually enhancing the degradation profile, which represents how the battery capacity declines over time, effectively modeling higher State of Health (SOH) levels. The results, shown in Figure~\ref{fig:deg}, indicate that if the initial SOH of the SL BESS exceeds approximately 0.78, and the degradation profile is improved by 20\% compared to the one reported in Figure~\ref{fig:deg_bat}, the model begins to favor battery installation. Installed capacity increases sharply as SOH improves, eventually reaching 0.8~MWh. Therefore, while cost remains the dominant factor, degradation performance is a decisive constraint that can balance feasibility and infeasibility.

Following this, a complementary sensitivity analysis was carried out to examine the impact of SL BESS cost on the installation decision. In this case, battery prices were systematically varied while keeping the electricity price profile and the degradation behavior fixed (assuming a 20\% improvement over current SL BESS performance). The results, shown in Figure~\ref{fig:cp}, show a sharp response of the model to price variations. For SL BESS prices below approximately 158~\euro{}/MWh, the system always installs the maximum storage capacity. As the price increases beyond that point, the installed capacity begins to decline progressively, reaching zero when the cost exceeds 170~\euro{}/MWh. The dashed vertical lines in the figure correspond to the SL BESS prices defined in scenarios C3 (green) and C2 (red), as reported in Table~\ref{table:battery_prices}. These findings confirm that, even under improved degradation conditions, cost remains a critical factor in determining the feasibility of battery deployment.

\subsection{Electricity market prices fluctuations}
\label{subsec:ele}

\begin{figure}[H]
    \centering
    \includegraphics[width=1\textwidth]{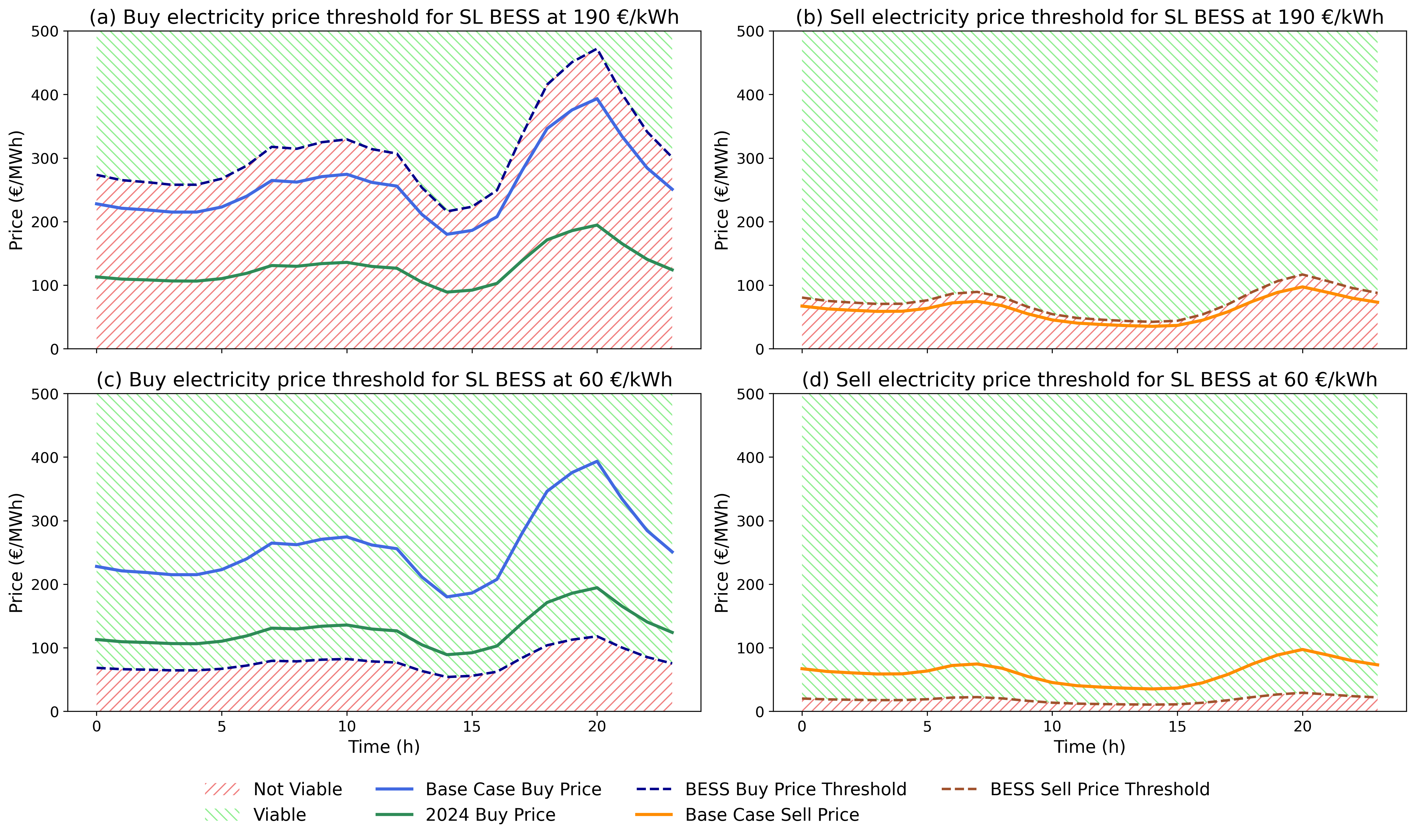} 
    \caption{Comparison of buy and sell threshold electricity prices for C2 (top row) and C4 (bottom row).}
    \label{fig:thresholds_combined}
\end{figure}

Figure~\ref{fig:thresholds_combined} illustrates the threshold electricity prices required to make the installation of a second-life (SL) BESS economically viable, based on hourly average prices. The left column shows the electricity buy price curves, while the right column displays the corresponding sell price curves. Each subplot includes three lines: the average electricity price for 2024, the base case price profile, and the dashed threshold curve that defines the minimum (or maximum) price required for battery investment to be profitable. The top row corresponds to a commercial SL BESS cost of 190~\euro{}/MWh (scenario C2), while the bottom row assumes a substantially reduced cost of 60~\euro{}/MWh (scenario C4).

As observed in the top row (C2), electricity prices would need to exceed those of the 2022 base case, already considered high, to make SL BESS investment profitable. This indicates that, at current commercial prices, and under assumptions regarding initial SoH 70\% and degradation behavior, second-life batteries remain economically unviable. In contrast, when SL BESS prices are reduced to 60~\euro{}/MWh (bottom row), the viability region expands considerably. Even under relatively low electricity market prices (e.g., similar to or below those of 2024), community-scale battery investments become attractive. This analysis reinforces that cost reductions in SL battery technologies are essential to unlocking their widespread deployment in community energy systems.


\subsection{Penetration of DER and demand}
\label{subsec:pen/dem}



\begin{figure}[H]
    \begin{minipage}[b]{0.48\textwidth}
        \centering
        \includegraphics[width=\textwidth]{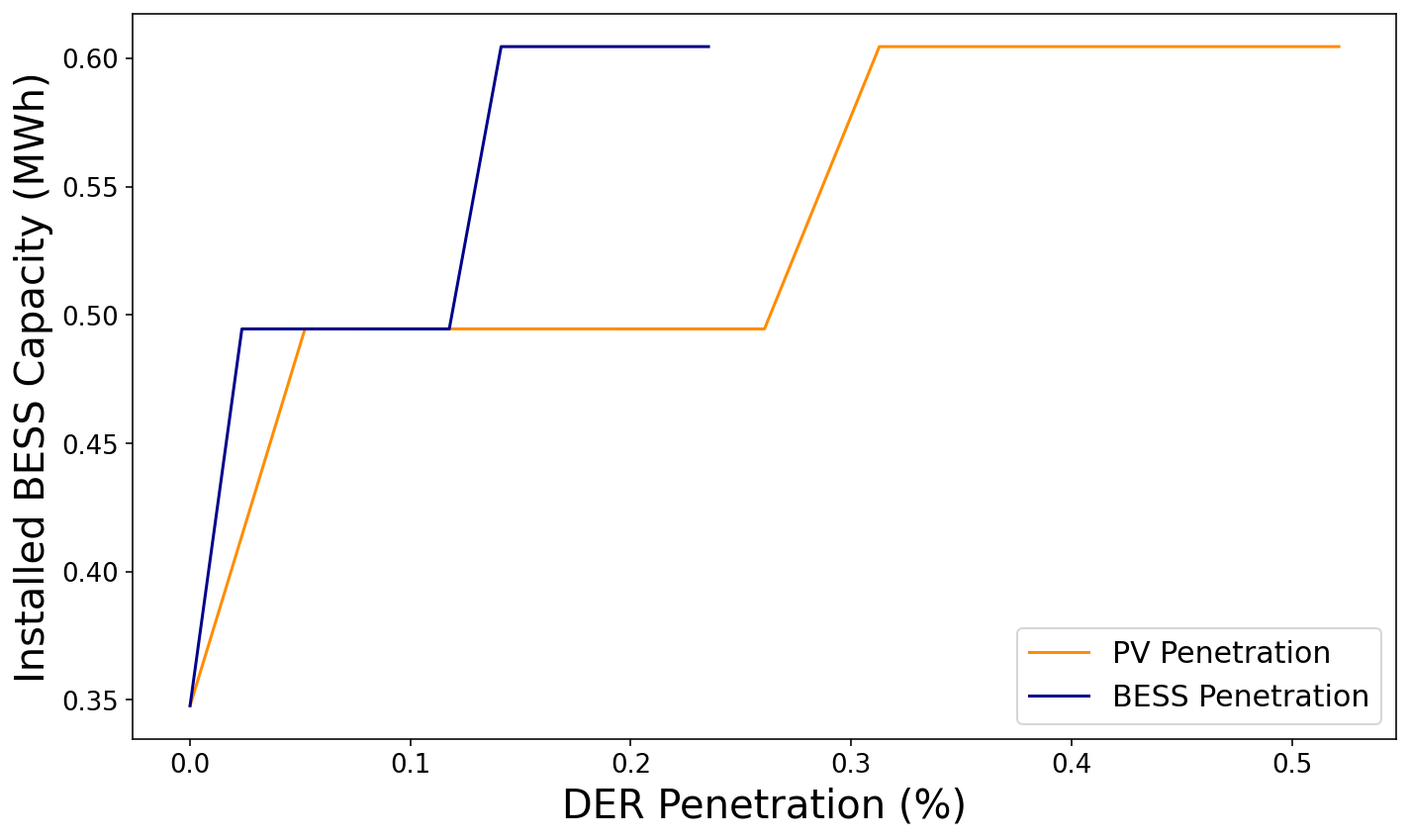}
        \caption{Installed BESS Capacity vs DER penetration}
        \label{fig:pen}
    \end{minipage}
    \hfill
    \begin{minipage}[b]{0.48\textwidth}
        \centering
        \includegraphics[width=\textwidth]{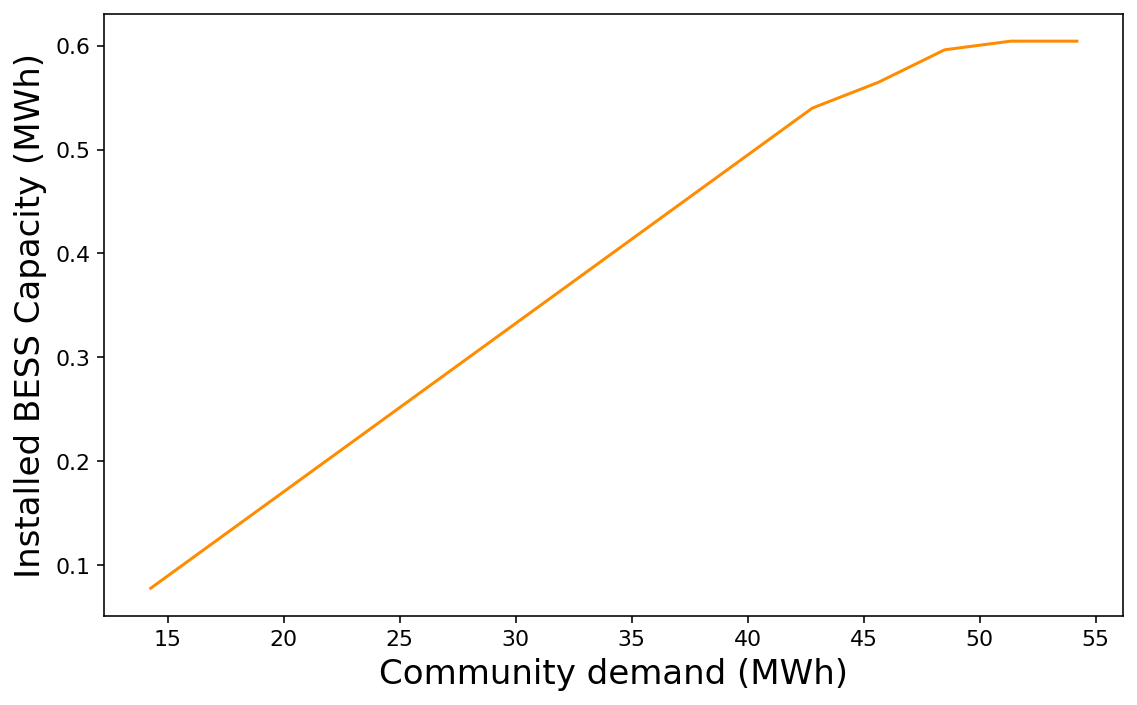}
        \caption{Installed BESS Capacity vs Community demand}
        \label{fig:dem}
    \end{minipage}
\end{figure}

This section investigates how the level of DER penetration and the total community demand influence the sizing of a shared SL BESS. Battery investment costs correspond to C3, which assumes a commercial price of 190~\euro{}/MWh for SL technologies. To analyze DER penetration, a continuous scaling factor is applied to reduce or increase the installed capacity of user-owned PV and BESS. As shown in Figure~\ref{fig:pen}, DER penetration is defined as the ratio between the installed capacity of individual DERs (PV and BESS) and the peak electricity demand of the community. A value of 0 indicates no user owns DERs, while 0.5 corresponds to the maximum penetration level considered in this study. The results show a clear trend: as DER penetration increases, the installed capacity of the shared community BESS also increases. Notably, the installation becomes economically viable even at moderate DER penetration levels of BESS. This suggests that the progressive deployment of distributed resources, both PV and BESS, could creates operational conditions that enhance the value of a centralized, shared battery system.

Regarding demand, Figure~\ref{fig:dem} illustrates the sensitivity of installed BESS capacity to changes in total community demand. The analysis shows that as demand increases, the incentive to install a community-scale battery also grows. The relationship is steeper at lower demand levels, meaning that initial increases in demand significantly influence battery sizing. However, this effect gradually plateaus, indicating that once demand surpasses a certain threshold, additional increments have a diminishing impact on the optimal battery capacity. This behavior suggests that community-scale storage becomes progressively more justifiable as demand scales up, but eventually stabilizes once basic flexibility and arbitrage needs are met.

\subsection{Market price fluctuations peak/off-peak hours}\label{subsec:peak/valle}

\begin{figure}[H]
    \centering
    \includegraphics[width=0.55\textwidth]{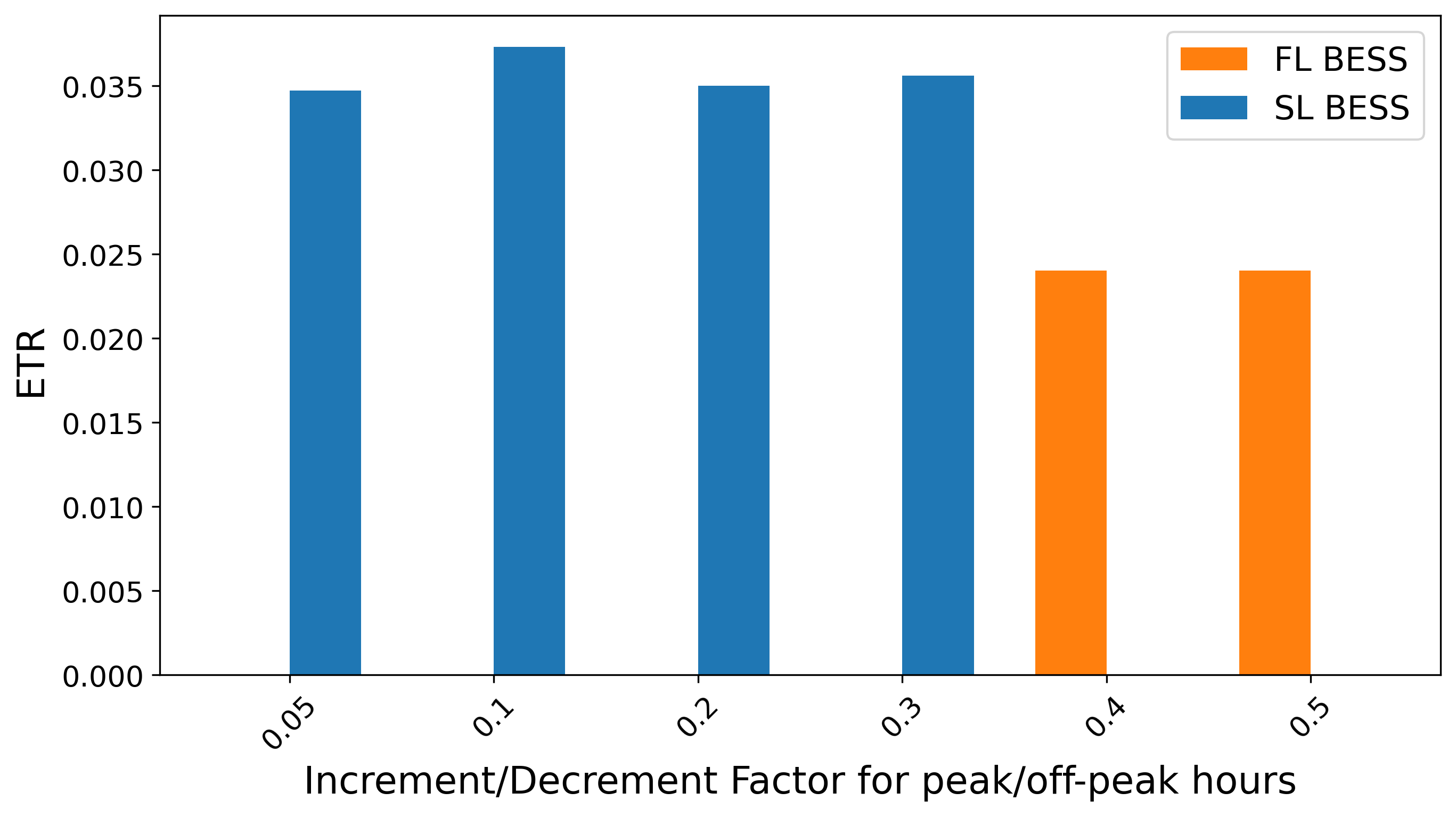}
    \caption{Installed BESS capacity and technology with modifications in peak/off-peak prices. ETR: Energy Traded-to-Demand Ratio}
    \label{fig:vp}
\end{figure}

This section analyzes how the installed capacity and selected technology of the community battery change as the electricity price profile becomes more pronounced—i.e., as peak prices rise and off-peak prices fall. To capture the effect on P2P trading, we also compute the Energy Traded-to-Demand Ratio (ETR), which is defined as the share of electricity exchanged among users relative to the total community demand. Across all cases, this ratio remains relatively low, averaging around 3.5\%, indicating that only a small portion of the total consumption is met through local energy exchanges. Thus, as the price profile becomes more accentuated, particularly beyond an increment/decrement factor of 0.3, the model shifts from selecting SL BESS to FL BESS technologies, as shown in Figure~\ref{fig:vp}. This transition suggests that under steeper price differentials, battery degradation becomes more influential than investment cost in the decision-making process. Consequently, the ETR drops from 3.5\% to 2.5\%, reflecting a reduced reliance on P2P trading due to the model’s preference for more efficient but costlier storage technologies. While this may appear counterintuitive, it reflects a shift in operational strategy rather than a reduction in system performance. In essence, the battery is used more selectively, prioritizing higher-value exchanges.




\begin{figure}[H]
    \centering
    \includegraphics[width=0.55\textwidth]{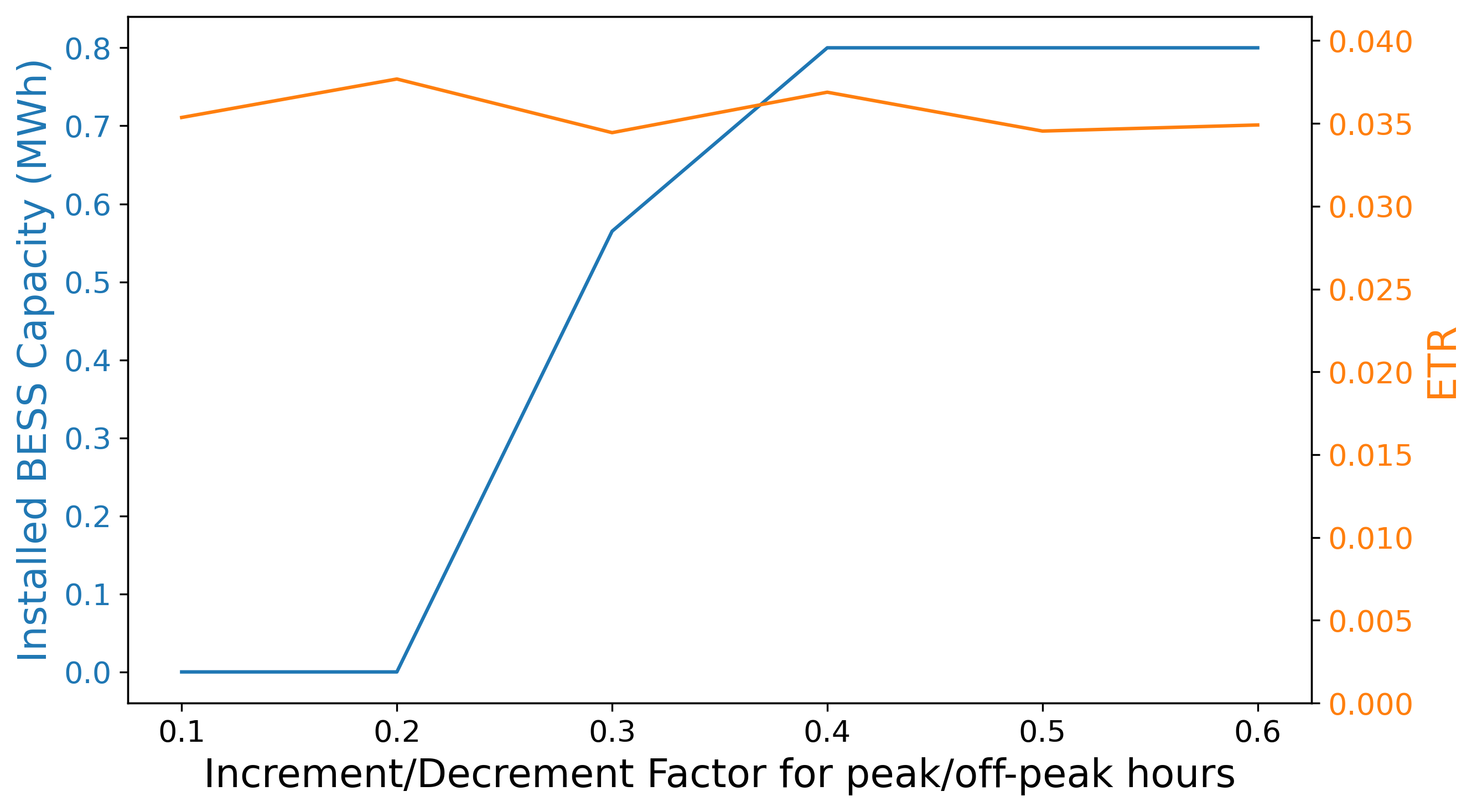}
    \caption{Installed BESS capacity with modifications in peak/off-peak prices}
    \label{fig:v}
\end{figure}

As previously explained, SL BESS prices from C3 in Table~\ref{table:battery_prices} led to the installation of the full storage capacity. For this reason, the prices from C2 are used in another analysis to capture the sensitivity of investment decisions, which are represented in \ref{fig:v}. The results indicate that investing in battery storage becomes progressively more attractive as the spread between peak and off-peak electricity prices increases. Plus, the system allocates capacity to storage once the price differential surpasses a certain threshold. This behavior aligns with economic logic, as larger spreads offer more profitable opportunities for charge and discharge cycles. However, despite the increase in installed capacity, the ETR remained relatively stable across scenarios. Therefore, while higher price volatility encourages battery investment, it does not necessarily lead to significantly greater energy trading.

\subsection{Energy trading and User BESS technology}
\label{subsec:etr/bat}



To assess how the characteristics of user-owned batteries influence the deployment of a community BESS and the extent of peer-to-peer energy exchanges, several representative scenarios were simulated. In each case, the Energy Trading Ratio (ETR) was computed as a proxy for the intensity of local energy trading. The results are summarized in Table~\ref{tab:trading}. A key insight is that the type and quality of the user-owned BESS significantly affect the trading dynamics. In Scenario 1, where users have FL batteries and the community installs a SL BESS, the ETR reaches 3.08\%, indicating a moderate level of exchange. However, in Scenario 2, where users rely on SL batteries with higher degradation, the ETR drops sharply to 0.52\%, even though a community BESS is still installed. This suggests that the degradation and reduced performance of individual batteries limit their contribution to local trading, thus reducing the availability of flexible resources. In contrast, Scenario 3 shows that when battery prices are significantly reduced (C4), making the FL technologies more competitive for community uses than the SL batteries. This scenario achieves a higher ETR of 3.57\%, emphasizing that improved efficiency and longevity of FL batteries enhance the potential for energy sharing. Furthermore, when the degradation profile of the SL BESS is artificially improved (Scenario 4), the ETR increases to 3.79\%, reinforcing that battery health and lifespan are central to sustaining local energy transactions.

Additionally, and based on the results presented in Table \ref{tab:trading}, the main conclusion is that SL batteries currently offer a competitive alternative to FL batteries, primarily due to their economic advantage, even assuming a 40\% price reduction. However, there are important limitations. In Scenario 3, where battery prices decrease significantly, investment cost is no longer a major constraint, and other factors such as efficiency become more decisive. In these cases, FL batteries are preferred. Moreover, in Figure \ref{fig:vp}, where the gap between peak and off-peak prices widened, battery degradation becomes more critical. Under such conditions, the model tends to favor first-life batteries due to their superior performance and longer lifespan.

\begin{table}[H]
\centering
\caption{Summary of configurations and resulting ETR in different scenarios}
\resizebox{\textwidth}{!}{%
\begin{tabular}{|c|c|c|c|c|c|c|}
\hline
\textbf{Sc.} & \textbf{BESS Price} & \textbf{Grid Energy Price} & \textbf{Community BESS} & \textbf{User BESS} & \textbf{BESS Tech.} & \textbf{ETR (\%)} \\
\hline
1 & C3 & BC & Estimated degradation & FL & SL & 3.08 \\
2 & C3 & BC & Estimated degradation & SL & SL & 0.52 \\
3 & C4 & BC & Estimated degradation & FL & FL (LFP/Gr) & 3.57 \\
4 & C3 & BC & Improved degradation & FL & SL & 3.79 \\
\hline
\end{tabular}%
}
\label{tab:trading}
\end{table}

\subsection{Computational performance analysis}
\label{subsec:comp}

This section analyzes the computational effort required to solve the optimization problem under different configurations. Table~\ref{tab:battery_time} summarizes the number of constraints, variables, and the average computation time as the number of battery technologies available for selection increases (from one to four). Each row corresponds to a scenario in which the model must choose among 1, 2, 3, or 4 battery technologies, respectively. The total number of variables and constraints is already substantial, exceeding 500,000 in all cases. This confirms the inherent complexity of the problem, particularly due to its hourly resolution, nonlinear degradation modeling, and the inclusion of investment and operational decisions.

Interestingly, the increase in the number of available batteries does not lead to a proportional rise in the number of constraints or variables. The increments are relatively marginal, indicating that the structure of the problem remains relatively stable despite the added decision layers. In terms of computational time, the results show no clear monotonic trend. For instance, when the model considers up to three battery options, the average solution time peaks at 94.64 seconds, suggesting that the solver faces greater difficulty identifying the optimal trade-off among these three configurations. However, when a fourth battery is introduced, the computation time decreases to 69.43 seconds. This nonlinearity suggests that solution time is influenced more by the problem landscape and the feasibility of configurations than by problem size alone. Some options may be quickly discarded due to high infeasibility, while others require deeper exploration. Thus, performance is shaped less by the number of decision variables and more by the complexity and distinguishability of the feasible solutions.

\begin{table}[H]
\centering
\begin{tabular}{|c|c|c|c|}
\hline
\textbf{Nº BESS} & \textbf{Constraints} & \textbf{Variables} & \textbf{Average Time (s)} \\
\hline
1 & 508,563 & 528,723 & 57.69 \\
2 & 510,004 & 529,685 & 72.27 \\
3 & 511,445 & 530,647 & 94.64 \\
4 & 512,886 & 531,609 & 69.43 \\
\hline
\end{tabular}
\caption{Average computation time as the number of batteries increases.}
\label{tab:battery_time}
\end{table}

\section{Conclusions}\label{sec:conclusions}
This study proposed a MINLP model to determine the optimal sizing of shared PV and BESS installations in an energy community where users already own individual PV and/or BESS systems and engage in P2P energy trading. The model accounts for multiple battery technologies, including FL and SL options, explicitly modeling degradation effects over time. The model was tested on the LV-DN with 206 buses.

The results show that, under current conditions, namely 2023–2024 electricity prices and commercial battery costs, the installation of community PV systems is economically viable, while deploying community BESS remains unprofitable, even for SL batteries, due to their still-elevated capital costs and degradation profile. Specifically, for SL batteries, economic viability requires the initial SoH to exceed approximately 0.78 and battery prices to remain below 160 €/MWh. Likewise, the conditions under which BESS deployment becomes economically justified, several sensitivity analyses were conducted. The results reveal that the most influential factors are: (1) the market price of BESS technologies, (2) electricity market prices, particularly purchase tariffs, and (3) the degradation profile of the BESS technology. Secondary factors, such as the penetration level of DERs among users and the aggregated peak demand of the community, also play a role but have comparatively less impact.

The analysis further highlights that SL batteries may become a competitive alternative to FL technologies if their degradation performance improves and/or their investment cost decreases significantly. Moreover, scenarios with high volatility in electricity prices, characterized by wider peak/off-peak gaps, favor using FL batteries, as degradation becomes more detrimental to economic performance than upfront costs. Finally, future work will incorporate uncertainties in generation, demand, and market prices to enhance the robustness and applicability of the proposed planning framework under real-world conditions.






\section*{Acknowledgment}
This work has been supported by ANID FONDECYT Iniciación 11240745.

\newpage
\bibliographystyle{elsarticle-num-names}
\bibliography{references} 

\section{Appendix} 

\subsection*{Degradation Model for LFP/Gr and LMO/Gr}

\begin{equation}
k_{\text{cal}} = |p_1| \cdot \exp\left( \frac{p_2}{T_{\text{degK}}} \right) \cdot \exp\left( \frac{p_3 \cdot U_a}{T_{\text{degK}}} \right)
\end{equation}

\begin{equation}
k_{\text{cyc}} = (p_4 + p_5 \cdot \text{DOD} + p_6 \cdot C_{\text{rate}}) \cdot \left[ \exp\left( \frac{p_7}{T_{\text{degK}}} \right) + \exp\left( \frac{-p_8}{T_{\text{degK}}} \right) \right]
\end{equation}

\begin{align}
\Delta q_{\text{t}} &= \text{scale} \cdot \left( q_{\text{t,prev}} + k_{\text{cal}} \cdot (\Delta t)^{p_{\text{cal}}} \right) \\
\Delta q_{\text{EFC}} &= \text{scale} \cdot \left( q_{\text{EFC,prev}} + k_{\text{cyc}} \cdot (\Delta \text{EFC})^{p_{\text{cyc}}} \right)
\end{align}

\begin{equation}
q = 1 - q_{\text{t}} - q_{\text{EFC}}
\end{equation}

\begin{table}[h!]
    \centering
    \caption{Model parameters and stressors used in degradation calculation.}
    \begin{tabular}{|c|p{10cm}|}
        \hline
        \textbf{Symbol} & \textbf{Description} \\
        \hline
        $T_{\text{degK}}$ & Absolute temperature in Kelvin \\
        $U_a$ & Anode potential versus lithium \\
        $\text{DOD}$ & Depth-of-discharge \\
        $C_{\text{rate}}$ & Average C-rate (charge/discharge current normalized by capacity) \\
        $p_1$ to $p_3$ & Empirical coefficients for calendar ageing \\
        $p_4$ to $p_8$ & Empirical coefficients for cycle ageing \\
        $p_{\text{cal}}$ & Power-law exponent for calendar degradation accumulation \\
        $p_{\text{cyc}}$ & Power-law exponent for cycle degradation accumulation \\
        $\Delta t$ & Elapsed time in days \\
        $\Delta \text{EFC}$ & Elapsed equivalent full cycles \\
        $\text{scale}$ & Degradation scaling factor \\
        \hline
    \end{tabular}
    \label{tab:model_parameters}
\end{table}

\subsection*{Degradation Model for NMC/Gr}

\begin{equation}
k_{\text{cal}} = p_1 \cdot \exp\left( p_2 \cdot T_{\text{degK}}^3 \cdot \left( \frac{1}{U_a^{1/3}} \right) \right)
\end{equation}

\begin{equation}
k_{\text{cyc}} = \left| 
p_3 
+ p_4 \cdot T_{\text{degK}}^3 \cdot \sqrt{\text{DOD}} 
+ p_5 \cdot \exp\left( \frac{C_{\text{rate}}^2}{\sqrt{T_{\text{degK}}}} \right) 
+ p_6 \cdot \exp\left( \sqrt{\text{DOD}} \cdot T_{\text{degK}}^2 \cdot \sqrt{C_{\text{rate}}} \right) 
\right|
\end{equation}

\begin{equation}
\Delta q_t = \text{scale} \cdot \left( q_{\text{Loss}t,\text{prev}} + k{\text{cal}} \cdot \left( \frac{\Delta t_{\text{days}}}{1000} \right)^{p_{\text{cal}}} \right)
\end{equation}

\begin{equation}
\Delta q_{\text{EFC}} = \text{scale} \cdot \left( q_{\text{Loss}{\text{EFC}},\text{prev}} + k{\text{cyc}} \cdot \left( \frac{\Delta \text{EFC}}{10^4} \right)^{p_{\text{cyc}}} \right)
\end{equation}

\begin{equation}
q = 1 - q_{\text{Loss}t} - q{\text{Loss}_{\text{EFC}}}
\end{equation}

\begin{table}[h!]
\centering
\caption{Parameter definitions for the NMC/Gr degradation model~\cite{Gasper2023}}
\begin{tabular}{|l|p{11cm}|}
\hline
\textbf{Parameter} & \textbf{Description} \\
\hline
$T_{\text{degK}}$ & Temperature in Kelvin \\
$U_a$ & Anode potential vs. Li/Li\textsuperscript{+} \\
$\text{DOD}$ & Depth of discharge \\
$C_{\text{rate}}$ & Average C-rate of charge/discharge \\
$p_1$ & Pre-exponential factor for calendar ageing \\
$p_2$ & Temperature and anode potential dependency for calendar ageing \\
$p_{\text{cal}}$ & Power-law exponent for calendar ageing time dependence \\
$p_3$ & Constant term in cycle ageing degradation rate \\
$p_4$ & Temperature and DOD coupling term \\
$p_5$ & Exponential dependence on C-rate and inverse square root of temperature \\
$p_6$ & Exponential term coupling DOD, temperature, and C-rate \\
$p_{\text{cyc}}$ & Power-law exponent for cycle ageing \\
\hline
\end{tabular}
\label{tab:nmcgr_parameters}
\end{table}

\subsection*{Degradation Model for NMC/LTO}
\begin{equation}
\alpha = \alpha_0 \cdot \exp\left( \frac{\alpha_1}{T_{\text{degK}}} \right) \cdot \exp\left( \frac{\alpha_2 \cdot \text{SOC}}{T_{\text{degK}}} \right)
\end{equation}

\begin{equation}
\beta = \beta_0 \cdot \exp\left( \frac{\beta_1}{T_{\text{degK}}} \right) \cdot \exp\left( \frac{\beta_2 \cdot \text{SOC}}{T_{\text{degK}}} \right)
\end{equation}

\begin{equation}
\gamma = \left( \gamma_0 + \gamma_1 \cdot C_{\text{rate}} + \gamma_2 \cdot \text{DOD}^3 \right) \cdot \exp\left( \frac{\gamma_3}{T_{\text{degK}}} \right)
\end{equation}

\begin{align}
\Delta q_{\text{gain},t} &= \text{scale} \cdot \left( q_{\text{Gain}_t,\text{prev}} + \alpha \cdot \Delta t^{\alpha_p} \right) \\
\Delta q_{\text{loss},t} &= \text{scale} \cdot \left( q_{\text{Loss}_t,\text{prev}} + \beta \cdot \Delta t^{\beta_p} \right) \\
\Delta q_{\text{EFC}} &= \text{scale} \cdot \left( q_{\text{Loss}_{\text{EFC}},\text{prev}} + \gamma \cdot \Delta \text{EFC}^{\gamma_p} \right)
\end{align}

\begin{equation}
q = 1 - q_{\text{Loss}t} + q{\text{Gain}t} - q{\text{Loss}_{\text{EFC}}}
\end{equation}

\begin{table}[h!]
\centering
\caption{Parameter definitions for the NMC-LTO 10Ah degradation model}
\begin{tabular}{|l|p{11cm}|}
\hline
\textbf{Parameter} & \textbf{Description} \\
\hline
$T_{\text{degK}}$ & Temperature in Kelvin \\
SOC & Average state of charge \\
DOD & Depth of discharge \\
$C_{\text{rate}}$ & C-rate during cycling \\
$\alpha_0$, $\alpha_1$, $\alpha_2$ & Coefficients for calendar ageing gain \\
$\alpha_p$ & Power-law exponent for calendar ageing gain \\
$\beta_0$, $\beta_1$, $\beta_2$ & Coefficients for calendar ageing loss \\
$\beta_p$ & Power-law exponent for calendar ageing loss \\
$\gamma_0$, $\gamma_1$, $\gamma_2$, $\gamma_3$ & Coefficients for cycle ageing degradation \\
$\gamma_p$ & Power-law exponent for cycle ageing degradation \\
\hline
\end{tabular}
\label{tab:nmc_lto_parameters}
\end{table}

\end{document}